\documentclass[reqno,b5paper]{amsart}

\usepackage{amssymb}
\usepackage{amsthm}
\usepackage{amsmath}
\usepackage{amsfonts}
\usepackage{mathrsfs}
\usepackage{enumerate}
\usepackage[mathscr]{eucal}

\DeclareMathOperator{\RE}{Re}
\DeclareMathOperator{\IM}{Im} 
 
 \DeclareMathOperator{\li}{li}

 \DeclareMathOperator{\Ad}{Ad}

\DeclareMathOperator{\II}{i}

\DeclareMathOperator{\RES}{Res}

\newtheorem{theorem}{Theorem}

\numberwithin{equation}{section}
\numberwithin{theorem}{section}

\newcommand{\set}[1]{\left\{#1\right\}}
\newcommand{\abs}[1]{\left\vert#1\right\vert}
\newcommand{\br}[1]{\left(#1\right)}
\newcommand{\SqBr}[1]{\left[#1\right]}
\newcommand{\spbr}[1]{\left(#1\right]}
\newcommand{\sppbr}[1]{\left[#1\right)}

\begin{document}

\title[On refinements of Gallagherian PGT's]{On refinements of rank one Gallagherian prime geodesic theorems}
\author[D\v zenan Gu\v si\'c]{D\v zenan Gu\v si\'c}

\newcommand{\acr}{\newline\indent}

\address{Department of Mathematics, University of Sarajevo, Zmaja od Bosne 35,\acr 71000 Sarajevo, Bosnia and Herzegovina}

\email{dzenang@pmf.unsa.ba}


\begin{abstract}
In his recent research, the author improved the error term in the prime geodesic theorem for compact, even-dimensional, rank one locally symmetric spaces. It turned out that the obtained estimate\\ $O\br{x^{2\rho-\frac{\rho}{n}}\br{\log x}^{-1}}$ coincides with the best known results for compact Riemann surfaces, three manifolds, and manifolds with cusps, where $n$ stands for the dimension of the space, and $\rho$ is the half-sum of positive roots. The above bound was then reduced to\\ $O\br{x^{2\rho-\rho\frac{2\cdot\br{2n}+1}{2n\cdot\br{2n}+1}}\br{\log x}^{\frac{n-1}{2n\cdot\br{2n}+1}-1}\br{\log\log x}^{\frac{n-1}{2n\cdot\br{2n}+1}+\varepsilon}}$ in the Gallagherian sense, with $\varepsilon$ $>$ $0$, and the key role played by the counting function $\psi_{2n}\br{x}$. The purpose of this research is to prove that the latter $O$-term can be further reduced. To do so, we derive new explicit formulas for the functions $\psi_{j}\br{x}$, $j$ $\geq$ $n$, and conditional formula for $\psi_{n-1}\br{x}$. Applying the Gallagher-Koyama techniques, we deduce the asymptotics for $\psi_{0}\br{x}$, and the Gallagherian prime geodesic theorems. The obtained error terms $O\br{x^{2\rho-\rho\frac{2j+1}{2nj+1}}\br{\log x}^{\frac{n-1}{2nj+1}-1}\br{\log\log x}^{\frac{n-1}{2nj+1}+\varepsilon}}$, $n-1$ $\leq$ $j$ $<$ $2n$, improve the $O$-term given above, with the optimal unconditional (conditional) size achieved for $j$ $=$ $n$ ($j$ $=$ $n-1$). If $j$ $=$ $n$ $\geq$ $4$, our new bound coincides with the best known estimate in the manifolds with cusps case. If $j$ $=$ $n-1$, the $O$-term fully agrees with the results in the Riemann surface case ($n$ $=$ $2$, $\rho$ $=$ $\frac{1}{2}\br{n-1}$ $=$ $\frac{1}{2}$), and the three manifolds case ($n$ $=$ $2$, $\rho$ $=$ $1$). Finally, for $j$ $=$ $n-1$, $n$ $\geq$ $4$, our result improves the best known bound in the manifolds with cusps case.
\end{abstract}

\maketitle

\section{Introduction}
As it is well known, there is a lot of literature on Selberg zeta functions. Most of papers, starting with the Selberg's one \cite{73Bunke}, discuss the Riemann surface case (two-dimensional case), or at most the finite volume case.

Quite opposite, only a few of papers treat more general settings: Gangolli \cite{24Bunke} (higher dimensional spherical case), Scott \cite{70Bunke} (three-dimensional hyperbolic manifolds case), Wakayama \cite{82Bunke} (general, locally symmetric spaces case), Gon-Park \cite{GonPark} (hyperbolic manifolds with cusps case), Gangolli-Warner \cite{25Bunke} (spherical rank one case), Moscovici-Stanton \cite{50Bunke} (locally symmetric manifolds case), Deitmar \cite{16Bunke} (higher rank locally symmetric spaces case). There are also papers like: Williams \cite{92Bunke}, Schuster \cite{69Bunke}, Parnovskij \cite{54Bunke}, etc.

For the followers of the Selberg's original work, who first obtained a meromorphic continuation of the logarithmic derivative of the Selberg zeta function, it was only possible to show that a certain power of the zeta function is meromorphic. The problem was circumvented by Fried \cite{22Bunke}, who established the connection with hyperbolic dynamics, and expressed the Selberg zeta function in terms of the Ruelle zeta function \cite{64Bunke}, which was known to be a meromorphic one.

Following the approach via the trace formula, Bunke and Olbrich \cite{Bunke} investigated the Selberg zeta function associated with a locally homogeneous vector bundle over the unit sphere bundle of a compact locally symmetric space of rank one. They associated two eliptic operators, one over the locally symmetric space and one over its compact dual symmetric space, and identified the spectrum of both operators together with the set of singularities of the zeta function. The main part of the Chapter 3 in \cite{Bunke} is devoted to the analytical continuation of the Selberg zeta function to all of $\mathbb{C}$, the description of its singularities, its functional equation, and its representation via regularized determinants. The meromorphic continuation of the Selberg zeta function (and not of an integer power) is obtained by showing that all residues of the logarithmic derivative are integers (see, \cite[p. 112, Prop. 3.14]{Bunke}).

Regarding the Ruelle zeta function, the authors in \cite{Bunke} concluded the Subsection 3.1 with the observation due to Fried \cite{22Bunke} that the Ruelle zeta function is a product of the Selberg zeta functions (see also, \cite{21Bunke} and \cite{23Bunke}).

The key role in this research will be played by these zeta functions. The necessary notation and normalization will also follow from Bunke-Olbrich's work.

Let $Y$ $=$ $\Gamma$ $\backslash$ $G$ $/$ $K$ be a compact, $n$-dimensional ($n$ even), locally symmetric Riemannian manifold with strictly negative sectional curvature, where $G$ is a connected semi-simple Lie group of real rank one, $K$ is a maximal compact subgroup of $G$, and $\Gamma$ is a discrete co-compact torsion-free subgroup of $G$.

Through the rest of the paper, we shall assume that the Riemannian metric over $Y$ (induced from the Killing form) is normalized such that the sectional curvature of $Y$ varies between $-4$ and $-1$.

Denote by $X$ the universal covering of $Y$.

As it is known, $X$ is a Riemannian symmetric space of rank one, so it is either a real $H\mathbb{R}^{k}$ or a complex $H\mathbb{C}^{m}$, or a quaternionic hyperbolic space $H\mathbb{H}^{m}$, or the hyperbolic Cayley plane $H\mathbb{C}a^{2}$.

We have, $n$ $=$ $k$, $2m$, $4m$ and $16$, respectively.

Put $\rho$ $=$ $\frac{1}{2}\sum\limits_{\alpha\in\Phi^{+}\br{\mathfrak{g},\mathfrak{a}}}\dim\br{\mathfrak{n}_{\alpha}}\alpha$ to be the half-sum of the positive roots, where $\Phi^{+}\br{\mathfrak{g},\mathfrak{a}}$ is a system of positive roots, $\mathfrak{n}$ $=$ $\sum\limits_{\alpha\in\Phi^{+}\br{\mathfrak{g},\mathfrak{a}}}\mathfrak{n}_{\alpha}$ is the sum of the root spaces, $\mathfrak{g}$ is the Lie algebra of $G$ with the Cartan decomposition $\mathfrak{g}$ $=$ $\mathfrak{k}$ $\oplus$ $\mathfrak{p}$, and $\mathfrak{a}$ is a maximal abelian subspace of $\mathfrak{p}$.

It follows that $\rho$ $=$ $\frac{1}{2}\br{k-1}$, $m$, $2m+1$ and $11$ when $n$ $=$ $k$, $2m$, $4m$ and $16$, respectively.

Let $\mathfrak{g}$ $=$ $\mathfrak{k}$ $\oplus$ $\mathfrak{a}$ $\oplus$ $\mathfrak{n}$ and $G$ $=$ $KAN$ be the Iwasawa decompositions.

Denote by $M$ the centralizer of $\mathfrak{a}$ in $K$ with the Lie algebra $\mathfrak{m}$.

The fact that $\Gamma$ is a co-compact and torsion-free subgroup of $G$ yields that there are only two types of conjugacy classes: the class of the identity $1\in\Gamma$, and classes of hyperbolic elements.

We put $C\Gamma$ to be the set of all conjugacy classes $\SqBr{\gamma}$ in $\Gamma$.

It is known that each hyperbolic $g\in G$ is conjugate to some $am$ $\in$ $A^{+}M$, where $A^{+}$ $=$ $\exp\br{\mathfrak{a}^{+}}$ for $\mathfrak{a}^{+}$ the positive Weyl chamber in $\mathfrak{a}$ (see, \cite{GangolliLength,25Bunke}).

By \cite[pp. 96-99]{Bunke}, the Selberg zeta function $Z_{S,\chi}\br{s,\sigma}$ is defined for $s$ $\in$ $\mathbb{C}$, $\RE\br{s}$ $>$ $\rho$ by the infinite product
\begin{align*}
&Z_{S,\chi}\br{s,\sigma}\\
=&\prod\limits_{\substack{1\neq\SqBr{g}\in C\Gamma\\\textrm{primitive}}}\prod\limits_{k=0}^{\infty}\det\br{1-\br{\sigma\br{m}\otimes\chi\br{g}\otimes S^{k}\br{\Ad\br{ma}_{\bar{\mathfrak{n}}}}}e^{-\br{s+\rho}l\br{g}}},
\end{align*}
where $S^{k}$ is the $k$-th symmetric power of an endomorphism, $\sigma$ and $\chi$ are finite-dimensional unitary representations of $M$ and $\Gamma$, $l\br{g}$ is the length of $g$, and $\bar{\mathfrak{n}}$ $=$ $\theta\mathfrak{n}$ with $\theta$ the Cartan involution of $\mathfrak{g}$.

For $s$ $\in$ $\mathbb{C}$, $\RE\br{s}$ $>$ $2\rho$, the Ruelle zeta function $Z_{R,\chi}\br{s,\sigma}$ is defined by the infinite product
\begin{align*}
Z_{R,\chi}\br{s,\sigma}=\prod\limits_{\substack{1\neq\SqBr{g}\in C\Gamma\\\textrm{primitive}}}\det\br{1-\br{\sigma\br{m}\otimes\chi\br{g}}}^{\br{-1}^{n-1}}.
\end{align*}

In view of Fried's observation noted above, we have that there are sets\\ $I_{p}$ $=$ $\set{\br{\tau,\lambda}\,:\,\tau\in\hat{M},\,\lambda\in\mathbb{R}}$, such that
\begin{align*}
Z_{R,\chi}\br{s,\sigma}=\prod\limits_{p=0}^{n-1}\prod\limits_{\br{\tau,\lambda}\in I_{p}}Z_{S,\chi}\br{s+\rho-\lambda,\tau\otimes\sigma}^{\br{-1}^{p}},
\end{align*}
where $\hat{M}$ (\,$\hat{\Gamma}$\,) stands for the unitary dual of $M$ (\,$\Gamma$\,).

Generally, by a prime geodesic theorem we mean an asymptotic formula for the lengths of the closed geodesics on a Riemannian manifold.

In his most recent research \cite{MDPI1}, the author improved the error term in DeGeorge's more than forty-year-old prime geodesic theorem \cite{DeGe}, up to\\ $O\br{x^{2\rho-\frac{\rho}{n}}\br{\log x}^{-1}}$, and then further reduced it in the Gallagherian sense \cite{Gal1,Gal2}, i.e., outside a set of finite logarithmic measure to\\ $O\br{x^{2\rho-\rho\frac{4n+1}{4n^{2}+1}}\br{\log x}^{\frac{n-1}{4n^{2}+1}-1}\br{\log\log x}^{\frac{n-1}{4n^{2}+1}+\varepsilon}}$, $\varepsilon$ $>$ $0$.

It was proven that there is a set $E$ of finite logarithmic measure, such that
\begin{align*}
\pi_{\Gamma}\br{x}=&\sum\limits_{p=0}^{n-1}\br{-1}^{p+1}\sum\limits_{\br{\tau,\lambda}\in I_{p}}\sum\limits_{\substack{\alpha\in S_{p,\tau,\lambda}^{\mathbb{R}}\\2\rho-\rho\frac{4n+1}{4n^{2}+1}<\alpha\leq 2\rho}}\li\br{x^{\alpha}}+\\
&O\br{x^{2\rho-\rho\frac{4n+1}{4n^{2}+1}}\br{\log x}^{\frac{n-1}{4n^{2}+1}-1}\br{\log\log x}^{\frac{n-1}{4n^{2}+1}+\varepsilon}}
\end{align*}
as $x$ $\rightarrow$ $\infty$, $x$ $\notin$ $E$, where $S_{p,\tau,\lambda}^{\mathbb{R}}$ is the set of real singularities of\\ $Z_{S,\chi}\br{s+\rho-\lambda,\tau\otimes\sigma}$, and $\pi_{\Gamma}\br{x}$ is the function counting prime geodesics on $Y$ of length not larger than $\log x$.

The result was obtained with the use of Chebyshev type counting function $\psi_{2n}\br{x}$ of order $2n$.

Writing the $O$-term in the form
\begin{align*}
O\br{x^{2\rho-\rho\frac{2\cdot\br{2n}+1}{2n\cdot\br{2n}+1}}\br{\log x}^{\frac{n-1}{2n\cdot\br{2n}+1}-1}\br{\log\log x}^{\frac{n-1}{2n\cdot\br{2n}+1}+\varepsilon}},
\end{align*}
one easily comes to the conjecture that the $O$-term actually depends on $\psi_{2n}\br{x}$, or more generally, on the selection of $\psi_{j}\br{x}$, $j$ $\in$ $\mathbb{N}$.

Finding possible answers to these questions, as well as a possible answer to the question of whether or not the above $O$-term could be further reduced by varying $\psi_{j}\br{x}$, $j$ $\in$ $\mathbb{N}$, represent two of the three main motives to conduct this research at all.

The third motive comes from the knowledge that Hejhal \cite{Hejhal1,Hejhal2}, made use of the function $\psi_{1}\br{x}$ $=$ $\psi_{n-1}\br{x}$ to derive the classical prime geodesic theorems for compact hyperbolic surfaces and generic hyperbolic surfaces of finite volume ($n$ $=$ $2$), and the fact that Park \cite{Park} did the same thing using the function $\psi_{d-1}\br{x}$ $=$ $\psi_{n-1}\br{x}$ in the real hyperbolic manifolds with cusps case ($n$ $=$ $d$).

In view of the motives, we derive new explicit formulas for $\psi_{j}\br{x}$, $j$ $\geq$ $n$ (Theorem 1), and conditional asymptotic formula for $\psi_{n-1}\br{x}$ (Theorem 4). Hence, we derive the corresponding Gallagherian estimates for $\psi_{0}$ (Theorems 2 and 5). From there, we obtain improved prime geodesic theorems (Theorems 3 and 6) in a standard way, making use of lines in \cite[p. 102]{Park}. The results, combined, state that for $j$ $\geq$ $n-1$ and $\varepsilon$ $>$ $0$, there is a set $E$ of finite logarithmic measure, such that
\begin{equation}\label{trokuttokut}
\begin{aligned}
\pi_{\Gamma}\br{x}=&\sum\limits_{p=0}^{n-1}\br{-1}^{p+1}\sum\limits_{\br{\tau,\lambda}\in I_{p}}\sum\limits_{\substack{\alpha\in S_{p,\tau,\lambda}^{\mathbb{R}}\\2\rho-\rho\frac{2j+1}{2nj+1}<\alpha\leq 2\rho}}\li\br{x^{\alpha}}+\\
&O\br{x^{2\rho-\rho\frac{2j+1}{2nj+1}}\br{\log x}^{\frac{n-1}{2nj+1}-1}\br{\log\log x}^{\frac{n-1}{2nj+1}+\varepsilon}}
\end{aligned}
\end{equation}
as $x$ $\rightarrow$ $\infty$, $x$ $\in$ $\mathbb{R}$ $\backslash$ $E$.

This statement confirms the correctness of our conjecture, and, as it will be explained in detail below, gives a positive answer to the second of our questions, justifying also the use of the function $\psi_{n-1}\br{x}$.

Note that the result (\ref{trokuttokut}) is fully in line with the best known results obtained in the hyperbolic, finite volume, Riemann surface with cusps case, hyperbolic three manifolds case, and higher dimensional, real hyperbolic manifolds with cusps case.

Indeed, for a non-compact, hyperbolic Riemann surface with cusps, Theorem 3.1 in \cite[p. 28]{AvdAnalMath}, states that for $\varepsilon$ $>$ $0$, there is a set $E$ of finite logarithmic measure, such that
\begin{align*}
\pi_{\Gamma}\br{x}=\sum\limits_{\frac{7}{10}<\alpha\leq 1}\li\br{x^{\alpha}}+O\br{x^{\frac{7}{10}}\br{\log x}^{\frac{1}{5}-1}\br{\log\log x}^{\frac{1}{5}+\varepsilon}}
\end{align*}
as $x$ $\rightarrow$ $\infty$, $x$ $\in$ $\mathbb{R}$ $\backslash$ $E$, where $\alpha$ is a zero of the attached Selberg zeta function.

Obviously, the $O$-term follows form our (\ref{trokuttokut}) for the selection $j$ $=$ $n-1$, $n$ $=$ $2$ and $\rho$ $=$ $\frac{1}{2}\br{n-1}$ $=$ $\frac{1}{2}$.

By Theorem 1.2 in \cite[p. 691]{Erata}, for a hyperbolic three manifold, and $\varepsilon$ $>$ $0$, there is a set $E$ of finite logarithmic measure, so that
\begin{align*}
\pi_{\Gamma}\br{x}=\li\br{x^{2}}+\sum\limits_{n=1}^{M}\li\br{x^{s_{n}}}+O\br{x^{\frac{21}{13}}\br{\log x}^{-\frac{11}{13}}\br{\log\log x}^{\frac{2}{13}+\varepsilon}}
\end{align*}
as $x$ $\rightarrow$ $\infty$, $x$ $\in$ $\mathbb{R}$ $\backslash$ $E$, where $s_{1}$, $s_{2}$, ..., $s_{M}$ are the real zeros of the associated Selberg zeta function contained in the interval $\br{1,2}$.

Clearly, the error term coincides with our (\ref{trokuttokut}) for $j$ $=$ $n-1$, $n$ $=$ $3$ and $\rho$ $=$ $\frac{1}{2}\br{n-1}$ $=$ $1$.

Moreover, for a $d$-dimensional manifold with cusps, and $\varepsilon$ $>$ $0$, Theorem 2 in \cite[p. 3021]{AvdZenan} asserts that there exists a set $E$ of finite logarithmic measure, such that
\begin{equation}\label{delta1}
\begin{aligned}
\pi_{\Gamma}\br{x}=&\sum\limits_{d-1-\frac{1}{2}\br{d-1}\frac{2d+1}{2d^{2}+1}<s_{n}\br{k}\leq d-1}\li\br{x^{s_{n}\br{k}}}+\\
&O\br{x^{d-1-\frac{1}{2}\br{d-1}\frac{2d+1}{2d^{2}+1}}\br{\log x}^{\frac{d-1}{2d^{2}+1}-1}\br{\log\log x}^{\frac{d-1}{2d^{2}+1}+\varepsilon}}
\end{aligned}
\end{equation}
as $x$ $\rightarrow$ $\infty$, $x$ $\in$ $\mathbb{R}$ $\backslash$ $E$, where $\br{s_{n}\br{k}-k}\br{2d_{0}-k-s_{n}\br{k}}$ is a small eigenvalue in $\SqBr{0,\frac{3}{4}d_{0}^{2}}$ of $\Delta_{k}$ on $\pi_{\sigma_{k},\lambda_{n}\br{k}}$ with $s_{n}\br{k}$ $=$ $d_{0}$ $+$ $\II{}\lambda_{n}\br{k}$ or $s_{n}\br{k}$ $=$ $d_{0}$ $-$ $\II{}\lambda_{n}\br{k}$ in $\spbr{\frac{3}{2}d_{0},2d_{0}}$, $\Delta_{k}$ is the Laplacian acting on the space of $k$-forms over the manifold, $d_{0}$ $=$ $\frac{d-1}{2}$, and $\pi_{\sigma_{k},\lambda_{n}\br{k}}$ is the principal series representation.

The $O$-term in this result also follows from (\ref{trokuttokut}) for $j$ $=$ $n$, $n$ $=$ $d$ and $\rho$ $=$ $\frac{1}{2}\br{n-1}$ $=$ $\frac{1}{2}\br{d-1}$. Note that our result (\ref{trokuttokut}) improves the $O$-term in (\ref{delta1}) even more for the choice $j$ $=$ $n-1$, $n$ $=$ $d$ and $\rho$ $=$ $\frac{1}{2}\br{n-1}$ $=$ $\frac{1}{2}\br{d-1}$.

The classical (unconditional) prime geodesic theorems have been considered by many authors over the years, not always for the same underlying locally symmetric space. The authors like: DeGeorge \cite{DeGe}, Gangolli \cite{GangolliLength}, Gangolli-Warner \cite{25Bunke}, Park \cite{Park}, Hejhal \cite{Hejhal1,Hejhal2}, Randol \cite{Randol}, Deitmar-Pavey \cite{SL4}, Deitmar \cite{20MDPI1,23MDPI1}, Sarnak \cite{21MDPI1}, Iwaniec \cite{12MDPI1}, Luo-Sarnak \cite{26MDPI1}, Cai \cite{27MDPI1}, Soundararajan-Young \cite{28MDPI1}, Friedman-Jorgenson-Kramer \cite{18MDPI1}, and many others.

We conclude this section with the observation that taking $\psi_{j}\br{x}$, $j$ $\geq$ $n$, does not yield a better result than the one obtained in \cite[p. 9, Th. 2]{MDPI1} (see, Section 7). The mentioned result, however, coincides with the best known results in the Riemann surface case \cite{Randol}, hyperbolic three manifolds case \cite{Erata}, and real hyperbolic manifolds with cusps case \cite{AvdZenan}.

\section{Preliminary material}
We introduce only the necessary notation based on \cite{Bunke} and \cite{MDPI1} (see also, \cite{JMAA,Slovaka}).

For the sake of simplicity, we write $\gamma$ for an element of $C\Gamma$, and $\gamma_{0}$ for a primitive element.

In particular, if $\gamma$ and $\gamma_{0}$ appear in the same formula, it is understood that $\gamma_{0}$ is the primitive element underlying $\gamma$.

We fix $\sigma$ $=$ $1$, $\chi$ $=$ $1$, and reduce the notation by omitting to write them in the sequel.

Put $\psi_{0}\br{x}$ to be the sum $\sum\limits_{1\neq\SqBr{\gamma}\in C\Gamma,\,N\br{\gamma}\leq x}\Lambda\br{\gamma}$, where $N\br{\gamma}$ $=$ $e^{l\br{\gamma}}$, $\Lambda\br{\gamma}$ $=$ $\log N\br{\gamma_{0}}$, and $\psi_{j}\br{x}$ $=$ $\int\limits_{0}^{x}\psi_{j-1}\br{t}dt$, $j$ $\in$ $\mathbb{N}$.

Let $\mathcal{T}$ be the set of all $\tau$ $\in$ $\hat{M}$ occuring in the representation
\begin{equation}\label{4zvijezde}
\begin{aligned}
Z_{R}\br{s}=\prod\limits_{p=0}^{n-1}\prod\limits_{\br{\tau,\lambda}\in I_{p}}Z_{S}\br{s+\rho-\lambda,\tau}^{\br{-1}^{p}}.
\end{aligned}
\end{equation}

\section{Results}
\vspace{5mm}

\subsection{Explicit formulas for $\psi_{j}\br{x}$, $j$ $\geq$ $n$}
In this section we derive new explicit formulas for the counting functions $\psi_{j}\br{x}$, $j$ $\geq$ $n$.

The obtained explicit formulas will be applied in the proof of refined Gallagherian prime geodesic theorem attached to locally symmetric spaces described in this research.

The main result of the section is the following theorem.
\begin{theorem}\label{Th1}
Let $Y$ be as above. If $j$ $\geq$ $n$, then
\begin{equation}\label{cetiri}
\begin{aligned}
\psi_{j}\br{x}=&\sum\limits_{p=0}^{n-1}\br{-1}^{p+1}\sum\limits_{\br{\tau,\lambda}\in I_{p}}
\sum\limits_{\rho<\alpha\leq 2\rho}\frac{x^{\alpha+j}}{\prod\limits_{k=0}^{j}\br{\alpha+k}}+\\
&\sum\limits_{p=0}^{n-1}\br{-1}^{p+1}\sum\limits_{\substack{\br{\tau,\lambda}\in I_{p}\\ \lambda=2\rho}}
\sum\limits_{\RE\br{\alpha}=\rho}\frac{x^{\alpha+j}}{\prod\limits_{k=0}^{j}\br{\alpha+k}},
\end{aligned}
\end{equation}
where $\alpha$ is a singularity of the Selberg zeta function $Z_{S}\br{s+\rho-\lambda,\tau}$.
\end{theorem}
\begin{proof}
Reasoning in the same way as in \cite[p. 98]{Park}, we obtain for $c$ $>$ $2\rho$ (see, Theorems A and B in \cite[pages 18 and 31]{Ingham}, and the relation (3.4) in \cite[p. 97]{Bunke} with $\sigma$ $=$ $1$, $\chi$ $=$ $1$, even $n$)
\begin{align*}
\psi_{j}\br{x}=\frac{1}{2\pi\II{}}\int\limits_{c-\II{}\infty}^{c+\II{}\infty}-\frac{Z_{R}^{'}\br{s}}{Z_{R}\br{s}}\frac{x^{s+j}}{s\br{s+1}...\br{s+j}}ds.
\end{align*}

Suppose that there is no pole of the integrand of $\psi_{j}\br{x}$ over boundary\\ of the closed domain $R\br{T}$ given by $R\br{T}$ $=$ $\set{s\in\mathbb{C}\,:\,\abs{s}\leq T,\,\RE\br{s}\leq\rho}$ $\cup$ $\set{s\in\mathbb{C}\,:\,\rho\leq\RE\br{s}\leq c,\,-\tilde{T}\leq\IM\br{s}\leq\tilde{T}}$, where $T$ $\gg$ $0$ and $\tilde{T}$ $=$ $\sqrt{T^{2}-\rho^{2}}$.

For a fixed $0$ $<$ $\varepsilon$ $<$ $c$ $-$ $\rho$, the Cauchy residue theorem gives us
\begin{equation}\label{pet}
\begin{aligned}
&\frac{1}{2\pi\II{}}\int\limits_{c-\II{}\tilde{T}}^{c+\II{}\tilde{T}}-\frac{Z_{R}^{'}\br{s}}{Z_{R}\br{s}}\frac{x^{s+j}}{\prod\limits_{k=0}^{j}\br{s+k}}ds=
\end{aligned}
\end{equation}
\begin{align*}
&\frac{-1}{2\pi\II{}}\br{\int\limits_{c+\II{}\tilde{T}}^{\rho+\varepsilon+\II{}\tilde{T}}+\int\limits_{\rho+\varepsilon+\II{}\tilde{T}}^{\rho+\II{}\tilde{T}}+\int\limits_{C_{T}}+
\int\limits_{\rho-\II{}\tilde{T}}^{\rho+\varepsilon-\II{}\tilde{T}}+
\int\limits_{\rho+\varepsilon-\II{}\tilde{T}}^{c-\II{}\tilde{T}}}\br{-\frac{Z_{R}^{'}\br{s}}{Z_{R}\br{s}}\frac{x^{s+j}}{\prod\limits_{k=0}^{j}\br{s+k}}}ds\\
&+\sum\limits_{\alpha\in R\br{T}}\RES_{s=\alpha}\br{-\frac{Z_{R}^{'}\br{s}}{Z_{R}\br{s}}\frac{x^{s+j}}{\prod\limits_{k=0}^{j}\br{s+k}}},
\end{align*}
where $C_{T}$ denotes the circular part of the boundary of $R\br{T}$ with the anti-clockwise orientation.

First, we estimate the integrals $\int\limits_{\rho+\varepsilon+\II{}\tilde{T}}^{\rho+\II{}\tilde{T}}\cdot ds$ and $\int\limits_{\rho-\II{}\tilde{T}}^{\rho+\varepsilon-\II{}\tilde{T}}\cdot ds$ on the right hand side of (\ref{pet}).

Denote by $N_{S,p}^{\tau,\lambda}\br{y}$ the number of singularities $\rho_{S,p}^{\tau,\lambda}$ $=$ $-\rho$ $+$ $\lambda$ $+$ $\II{}\gamma_{S,p}^{\tau,\lambda}$ of $Z_{S}\br{s+\rho-\lambda,\tau}$ on the interval $-\rho$ $+$ $\lambda$ $+$ $\II{}x$, $0$ $<$ $x$ $\leq$ $y$. By Theorem 3.15 in \cite[p. 113]{Bunke} and Theorem 9.1 in \cite[p. 89]{Duistermat}, $N_{S,p}^{\tau,\lambda}\br{y}$ $=$ $C_{1}y^{n}$ $+$ $O\br{y^{n-1}\br{\log x}^{-1}}$ for some explicitly known constant $C_{1}$. Since the factor $\br{\log x}^{-1}$ does not improve our result, we take
\begin{equation}\label{sest}
\begin{aligned}
N_{S,p}^{\tau,\lambda}\br{y}=C_{1}y^{n}+O\br{y^{n-1}}.
\end{aligned}
\end{equation}

Furthermore, by Theorem 4.1, claim ($b$) in \cite[p. 314]{Slovaka}, for $\delta$ $>$ $0$, and $t$ $\gg$ $0$, such that $\II{}t$ is not a singularity of $Z_{S}\br{s,\tau}$ for all $\tau$ $\in$ $\mathcal{T}$,
\begin{equation}\label{sedam}
\begin{aligned}
\frac{Z_{R}^{'}\br{s}}{Z_{R}\br{s}}=O\br{t^{n-1+\delta}}+\sum\limits_{p=0}^{n-1}\br{-1}^{p}\sum\limits_{\substack{\br{\tau,\lambda}\in I_{p}\\ \lambda=2\rho}}
\sum\limits_{\abs{t-\gamma_{S,p}^{\tau,\lambda}}\leq 1}\frac{1}{s-\rho_{S,p}^{\tau,\lambda}}
\end{aligned}
\end{equation}
for $s$ $=$ $\sigma_{1}$ $+$ $\II{}t$, $\rho$ $\leq$ $\sigma_{1}$ $<$ $\frac{1}{4}t$ $-$ $\rho$, and
\begin{equation}\label{osam}
\begin{aligned}
\frac{Z_{R}^{'}\br{s}}{Z_{R}\br{s}}=O\br{\frac{1}{\eta}t^{n-1+\delta}}
\end{aligned}
\end{equation}
for $s$ $=$ $\sigma_{1}$ $+$ $\II{}t$, $\rho$ $+$ $\eta$ $\leq$ $\sigma_{1}$ $<$ $\frac{1}{4}t$ $-$ $\rho$, where $0$ $<$ $\eta$ $\leq$ $2\rho$.

Note that Park derived a variant of the equation (\ref{sedam}) and a variant of the equation (\ref{osam}) in the case of real hyperbolic manifolds with cusps (see, equations (3.9) and (3.10) in \cite[p. 99]{Park}).

Now, reasoning in the same way as Park did, that is, applying the relations (\ref{sedam}) and (\ref{sest}), we obtain (Cf. \cite[p. 99, relation (3.11)]{Park} and \cite[p. 368, relation (2)]{Koreja})
\begin{equation}\label{devet}
\begin{aligned}
\frac{1}{2\pi\II{}}\int\limits_{\rho+\II{}\tilde{T}}^{\rho+\varepsilon+\II{}\tilde{T}}-\frac{Z_{R}^{'}\br{s}}{Z_{R}\br{s}}\frac{x^{s+j}}{\prod\limits_{k=0}^{j}\br{s+k}}ds
=O\br{x^{j+\rho+\varepsilon}T^{-j-2+n+\delta}}.
\end{aligned}
\end{equation}

The integral $\int\limits_{\rho-\II{}\tilde{T}}^{\rho+\varepsilon-\II{}\tilde{T}}\cdot ds$ can be treated in the same way, and gives us
\begin{equation}\label{deset}
\begin{aligned}
\frac{1}{2\pi\II{}}\int\limits_{\rho+\varepsilon-\II{}\tilde{T}}^{\rho-\II{}\tilde{T}}-\frac{Z_{R}^{'}\br{s}}{Z_{R}\br{s}}\frac{x^{s+j}}{\prod\limits_{k=0}^{j}\br{s+k}}ds
=O\br{x^{j+\rho+\varepsilon}T^{-j-2+n+\delta}}.
\end{aligned}
\end{equation}

In order to estimate the integrals $\int\limits_{c+\II{}\tilde{T}}^{\rho+\varepsilon+\II{}\tilde{T}}\cdot ds$ and $\int\limits_{\rho+\varepsilon-\II{}\tilde{T}}^{c-\II{}\tilde{T}}\cdot ds$, which appear in (\ref{pet}), we apply the asymptotics (\ref{osam}).

It follows immediately that
\begin{equation}\label{jedanaest}
\begin{aligned}
\frac{1}{2\pi\II{}}\int\limits_{\rho+\varepsilon+\II{}\tilde{T}}^{c+\II{}\tilde{T}}-\frac{Z_{R}^{'}\br{s}}{Z_{R}\br{s}}\frac{x^{s+j}}{\prod\limits_{k=0}^{j}\br{s+k}}ds
=O\br{\varepsilon^{-1}x^{c+j}T^{-j-2+n+\delta}}
\end{aligned}
\end{equation}
and
\begin{equation}\label{dvanaest}
\begin{aligned}
\frac{1}{2\pi\II{}}\int\limits_{c-\II{}\tilde{T}}^{\rho+\varepsilon-\II{}\tilde{T}}-\frac{Z_{R}^{'}\br{s}}{Z_{R}\br{s}}\frac{x^{s+j}}{\prod\limits_{k=0}^{j}\br{s+k}}ds
=O\br{\varepsilon^{-1}x^{c+j}T^{-j-2+n+\delta}}.
\end{aligned}
\end{equation}

The remaining integral on the right hand side of (\ref{pet}) is $\int\limits_{C_{T}}\cdot ds$.

Recall the following facts.

By Corollary 4.2 in \cite[p. 530]{JMAA}, a meromorphic extension over $\mathbb{C}$ of the Ruelle zeta function $Z_{R}\br{s}$ can be expressed as
\begin{equation}\label{trinaest}
\begin{aligned}
Z_{R}\br{s}=\frac{Z_{R}^{1}\br{s}}{Z_{R}^{2}\br{s}},
\end{aligned}
\end{equation}
where $Z_{R}^{1}\br{s}$ and $Z_{R}^{2}\br{s}$ are entire functions of order at most $n$ over $\mathbb{C}$ (Cf. \cite[p. 91, Th. 1.1]{Park}).

Moreover, according to Fried's known result (see, \cite[p. 509, Prop. 7]{22Bunke}), there is a constant $C$ $>$ $0$ such that for arbitrarily large choices of $r$
\begin{equation}\label{cetrnaest}
\begin{aligned}
\int\limits_{r}\abs{\frac{Z^{'}\br{s}}{Z\br{s}}}\abs{ds}\leq Cr^{n}\log r,
\end{aligned}
\end{equation}
where $Z\br{s}$ is a ratio of two nonzero entire functions of order not larger than $n$.

Following Park's reasoning in the case of real hyperbolic manifolds with cusps \cite[p. 99]{Park}, i.e., using the assertions (\ref{trinaest}) and (\ref{cetrnaest}), we estimate

\begin{equation}\label{petnaest}
\begin{aligned}
\frac{1}{2\pi\II{}}\int\limits_{C_{T}}-\frac{Z_{R}^{'}\br{s}}{Z_{R}\br{s}}\frac{x^{s+j}}{\prod\limits_{k=0}^{j}\br{s+k}}ds
=&O\br{x^{\rho+j}T^{-j-1}\int\limits_{C_{T}}\abs{\frac{Z_{R}^{'}\br{s}}{Z_{R}\br{s}}}\abs{ds}}\\
=&O\br{x^{\rho+j}T^{-j-1}\int\limits_{\abs{s}=T}\abs{\frac{Z_{R}^{'}\br{s}}{Z_{R}\br{s}}}\abs{ds}}\\
=&O\br{x^{\rho+j}T^{-j-1+n}\log T}.
\end{aligned}
\end{equation}

Finally, we analyze the sum given by residues over $R\br{T}$ in (\ref{pet}).

Bearing in mind the fact that the Ruelle zeta function $Z_{R}\br{s}$ is a product of the corresponding Selberg zeta functions (see, relation (\ref{4zvijezde})), and the fact that the singularities are precisely described by Theorem 3.15 in \cite[p. 113]{Bunke}, we deduce (Cf. \cite[p. 110]{Park}, \cite[p. 369]{Koreja} and \cite[p. 744, Th. 4.6]{GonPark} in the real hyperbolic manifolds with cusps case),
\begin{equation}\label{sesnaest}
\begin{aligned}
&\sum\limits_{\alpha\in R\br{T}}\RES_{s=\alpha}\br{-\frac{Z_{R}^{'}\br{s}}{Z_{R}\br{s}}\frac{x^{s+j}}{\prod\limits_{k=0}^{j}\br{s+k}}}\\
=&\sum\limits_{\rho<\alpha\leq 2\rho}\RES_{s=\alpha}\br{-\frac{Z_{R}^{'}\br{s}}{Z_{R}\br{s}}\frac{x^{s+j}}{\prod\limits_{k=0}^{j}\br{s+k}}}+
\sum\limits_{\substack{\alpha\in R\br{T}\\ \RE\br{\alpha}=\rho}}\RES_{s=\alpha}\br{.}\\
&+\sum\limits_{\substack{\alpha\in R\br{T}\\ \RE\br{\alpha}<\rho}}\RES_{s=\alpha}\br{.}\\
=&\sum\limits_{p=0}^{n-1}\br{-1}^{p+1}\sum\limits_{\br{\tau,\lambda}\in I_{p}}\sum\limits_{\rho<\alpha\leq 2\rho}
\RES_{s=\alpha}\br{\frac{Z_{S}^{'}\br{s+\rho-\lambda,\tau}}{Z_{S}\br{s+\rho-\lambda,\tau}}\frac{x^{s+j}}{\prod\limits_{k=0}^{j}\br{s+k}}}
\end{aligned}
\end{equation}
\begin{align*}
&+\sum\limits_{p=0}^{n-1}\br{-1}^{p+1}\sum\limits_{\substack{\br{\tau,\lambda}\in I_{p}\\ \lambda=2\rho}}\sum\limits_{\substack{\alpha\in R\br{T}\\ \RE\br{\alpha}=\rho}}
\RES_{s=\alpha}\br{\frac{Z_{S}^{'}\br{s+\rho-\lambda,\tau}}{Z_{S}\br{s+\rho-\lambda,\tau}}\frac{x^{s+j}}{\prod\limits_{k=0}^{j}\br{s+k}}}+
\end{align*}
\begin{align*}
&\sum\limits_{\substack{\alpha\in R\br{T}\\ \RE\br{\alpha}<\rho}}\RES_{s=\alpha}\br{-\frac{Z_{R}^{'}\br{s}}{Z_{R}\br{s}}\frac{x^{s+j}}{\prod\limits_{k=0}^{j}\br{s+k}}}.
\end{align*}

Consider the first sum on the right hand side of (\ref{sesnaest}).

For arbitrarily selected and then fixed $p$ $\in$ $\set{0,1,...,n-1}$ and $\br{\tau,\lambda}$ $\in$ $I_{p}$, let's take a closer look at $\RES_{s=\alpha}\br{\frac{Z_{S}^{'}\br{s+\rho-\lambda,\tau}}{Z_{S}\br{s+\rho-\lambda,\tau}}\frac{x^{s+j}}{\prod\limits_{k=0}^{j}\br{s+k}}}$, where $\rho$ $<$ $\alpha$ $\leq$ $2\rho$.

Clearly, all such $\alpha$'s are the poles of $\frac{Z_{S}^{'}\br{s+\rho-\lambda,\tau}}{Z_{S}\br{s+\rho-\lambda,\tau}}\frac{x^{s+j}}{\prod\limits_{k=0}^{j}\br{s+k}}$, i.e., they are the poles of $\frac{Z_{S}^{'}\br{s+\rho-\lambda,\tau}}{Z_{S}\br{s+\rho-\lambda,\tau}}$. This means that they are simple poles of $\frac{Z_{S}^{'}\br{s+\rho-\lambda,\tau}}{Z_{S}\br{s+\rho-\lambda,\tau}}$, that is, they are the singularities of $Z_{S}\br{s+\rho-\lambda,\tau}$.

Let $\alpha$ be a singularity of $Z_{S}\br{s+\rho-\lambda,\tau}$, with $\rho$ $<$ $\alpha$ $\leq$ $2\rho$.

Following Hejhal's reasoning in the compact Riemann surface case \cite[pp. 88-89]{Hejhal1}, we write
\begin{align*}
\frac{Z_{S}^{'}\br{s+\rho-\lambda,\tau}}{Z_{S}\br{s+\rho-\lambda,\tau}}=\frac{o_{\alpha}}{s-\alpha}\br{1+\sum\limits_{i=1}^{+\infty}a_{i}\br{\alpha}\br{s-\alpha}^{i}},
\end{align*}
where $o_{\alpha}$ is the order of $\alpha$, and $a_{i}\br{\alpha}$'s are the corresponding coefficients.

We obtain,
\begin{align*}
&\RES_{s=\alpha}\br{\frac{Z_{S}^{'}\br{s+\rho-\lambda,\tau}}{Z_{S}\br{s+\rho-\lambda,\tau}}\frac{x^{s+j}}{\prod\limits_{k=0}^{j}\br{s+k}}}\\
=&\lim\limits_{s\rightarrow\alpha}\br{s-\alpha}\frac{Z_{S}^{'}\br{s+\rho-\lambda,\tau}}{Z_{S}\br{s+\rho-\lambda,\tau}}\frac{x^{s+j}}{\prod\limits_{k=0}^{j}\br{s+k}}=
o_{\alpha}\frac{x^{\alpha+j}}{\prod\limits_{k=0}^{j}\br{\alpha+k}}.
\end{align*}

We carry out the same argumentation in the case of the second sum on the right-hand side of (\ref{sesnaest}).

Hence,
\begin{equation}\label{sedamnaest}
\begin{aligned}
&\sum\limits_{\alpha\in R\br{T}}\RES_{s=\alpha}\br{-\frac{Z_{R}^{'}\br{s}}{Z_{R}\br{s}}\frac{x^{s+j}}{\prod\limits_{k=0}^{j}\br{s+k}}}\\
=&\sum\limits_{p=0}^{n-1}\br{-1}^{p+1}\sum\limits_{\br{\tau,\lambda}\in I_{p}}\sum\limits_{\rho<\alpha\leq 2\rho}\frac{x^{\alpha+j}}{\prod\limits_{k=0}^{j}\br{\alpha+k}}+\\
&\sum\limits_{p=0}^{n-1}\br{-1}^{p+1}\sum\limits_{\substack{\br{\tau,\lambda}\in I_{p}\\ \lambda=2\rho}}\sum\limits_{\substack{\alpha\in R\br{T}\\ \RE\br{\alpha}=\rho}}\frac{x^{\alpha+j}}{\prod\limits_{k=0}^{j}\br{\alpha+k}}+\\
&\sum\limits_{\substack{\alpha\in R\br{T}\\ \RE\br{\alpha}<\rho}}\RES_{s=\alpha}\br{-\frac{Z_{R}^{'}\br{s}}{Z_{R}\br{s}}\frac{x^{s+j}}{\prod\limits_{k=0}^{j}\br{s+k}}},
\end{aligned}
\end{equation}
where, as usual, we omit to write $o_{\alpha}$'s for the sake of clarity.

In the same way Theorem 1.1 in \cite[p. 91]{Park} asserts that the Ruelle zeta function in the real hyperbolic manifolds with cusps case is a meromorphic function of order equal to the dimension of the underlying locally symmetric space, now, the equality (\ref{trinaest}) asserts that the Ruelle zeta function $Z_{R}\br{s}$ is a meromorphic function of order not larger than $n$.

To deal with the poles $\alpha$ $\in$ $R\br{T}$, $\RE\br{\alpha}$ $<$ $\rho$ of $-\frac{Z_{R}^{'}\br{s}}{Z_{R}\br{s}}\frac{x^{s+j}}{\prod\limits_{k=0}^{j}\br{s+k}}$ in (\ref{sedamnaest}), we proceed in exactly the same way as in \cite[p. 100]{Park}. Thus, we define $Z\br{s}$ to be a meromorphic function of order not larger than $n$ obtained by removing the singularities of $Z_{R}\br{s}$ in the right half plane $\RE\br{s}$ $\geq$ $\rho$ using the canonical product. It follows that
\begin{align*}
&\sum\limits_{\substack{\alpha\in R\br{T}\\ \RE\br{\alpha}<\rho}}\RES_{s=\alpha}\br{-\frac{Z_{R}^{'}\br{s}}{Z_{R}\br{s}}\frac{x^{s+j}}{\prod\limits_{k=0}^{j}\br{s+k}}}\\
=&\sum\limits_{\substack{\alpha\in R\br{T}\\ \RE\br{\alpha}<\rho}}\RES_{s=\alpha}\br{-\frac{Z^{'}\br{s}}{Z\br{s}}\frac{x^{s+j}}{\prod\limits_{k=0}^{j}\br{s+k}}}=
\frac{1}{2\pi\II{}}\int\limits_{C^{'}}-\frac{Z^{'}\br{s}}{Z\br{s}}\frac{x^{s+j}}{\prod\limits_{k=0}^{j}\br{s+k}}ds,
\end{align*}
where $C^{'}$ $=$ $C_{T}$ $\cup$ $\set{s\in\mathbb{C}\,:\,\RE\br{s}=\rho,\,-\tilde{T}\leq\IM\br{s}\leq\tilde{T}}$ is taken with the anti-clockwise orientation.

Suppose that
\begin{align*}
C^{''}=&\set{s\in\mathbb{C}\,:\,\abs{s}=T,\,\RE\br{s}\geq\rho}\cup\\
&\set{s\in\mathbb{C}\,:\,\RE\br{s}=\rho,\,-\tilde{T}\leq\IM\br{s}\leq\tilde{T}}\\
=&C_{T}^{c}\cup \set{s\in\mathbb{C}\,:\,\RE\br{s}=\rho,\,-\tilde{T}\leq\IM\br{s}\leq\tilde{T}}
\end{align*}
is also taken with the anti-clockwise orientation.

We may write
\begin{align*}
&\sum\limits_{\substack{\alpha\in R\br{T}\\ \RE\br{\alpha}<\rho}}\RES_{s=\alpha}\br{-\frac{Z_{R}^{'}\br{s}}{Z_{R}\br{s}}\frac{x^{s+j}}{\prod\limits_{k=0}^{j}\br{s+k}}}\\
=&\frac{1}{2\pi\II{}}\int\limits_{C^{'}}-\frac{Z^{'}\br{s}}{Z\br{s}}\frac{x^{s+j}}{\prod\limits_{k=0}^{j}\br{s+k}}ds+\frac{x^{\rho+j}}{2\pi}\cdot0\\
=&\frac{1}{2\pi\II{}}\int\limits_{C^{'}}-\frac{Z^{'}\br{s}}{Z\br{s}}\frac{x^{s+j}}{\prod\limits_{k=0}^{j}\br{s+k}}ds+
\frac{x^{\rho+j}}{2\pi}\int\limits_{C^{''}}-\frac{Z^{'}\br{s}}{Z\br{s}}\frac{ds}{\prod\limits_{k=0}^{j}\br{s+k}}.
\end{align*}

Applying the same reasoning as in the derivation of (\ref{petnaest}), we deduce
\begin{align*}
&\abs{\frac{1}{2\pi\II{}}\int\limits_{C^{'}}-\frac{Z^{'}\br{s}}{Z\br{s}}\frac{x^{s+j}}{\prod\limits_{k=0}^{j}\br{s+k}}ds+
\frac{x^{\rho+j}}{2\pi}\int\limits_{C^{''}}-\frac{Z^{'}\br{s}}{Z\br{s}}\frac{ds}{\prod\limits_{k=0}^{j}\br{s+k}}}\\
\leq&\frac{x^{\rho+j}}{2\pi}\int\limits_{C^{'}}\abs{\frac{Z^{'}\br{s}}{Z\br{s}}}\frac{\abs{ds}}{\prod\limits_{k=0}^{j}\abs{s+k}}+
\frac{x^{\rho+j}}{2\pi}\int\limits_{C^{''}}\abs{\frac{Z^{'}\br{s}}{Z\br{s}}}\frac{\abs{ds}}{\prod\limits_{k=0}^{j}\abs{s+k}}=
\end{align*}
\begin{align*}
&\frac{x^{\rho+j}}{2\pi}\int\limits_{C_{T}}\abs{\frac{Z^{'}\br{s}}{Z\br{s}}}\frac{\abs{ds}}{\prod\limits_{k=0}^{j}\abs{s+k}}+
\frac{x^{\rho+j}}{2\pi}\int\limits_{\rho-\II{}\tilde{T}}^{\rho+\II{}\tilde{T}}\abs{\frac{Z^{'}\br{s}}{Z\br{s}}}\frac{\abs{ds}}{\prod\limits_{k=0}^{j}\abs{s+k}}+\\
&\frac{x^{\rho+j}}{2\pi}\int\limits_{C_{T}^{c}}\abs{\frac{Z^{'}\br{s}}{Z\br{s}}}\frac{\abs{ds}}{\prod\limits_{k=0}^{j}\abs{s+k}}+
\frac{x^{\rho+j}}{2\pi}\int\limits_{\rho+\II{}\tilde{T}}^{\rho-\II{}\tilde{T}}\abs{\frac{Z^{'}\br{s}}{Z\br{s}}}\frac{\abs{ds}}{\prod\limits_{k=0}^{j}\abs{s+k}}\\
=&\frac{x^{\rho+j}}{2\pi}\int\limits_{\abs{s}=T}\abs{\frac{Z^{'}\br{s}}{Z\br{s}}}\frac{\abs{ds}}{\prod\limits_{k=0}^{j}\abs{s+k}}
\leq C_{1}x^{\rho+j}T^{-j-1}\int\limits_{\abs{s}=T}\abs{\frac{Z^{'}\br{s}}{Z\br{s}}}\abs{ds}\\
\leq&C_{2}x^{\rho+j}T^{-j-1+n}\log T.
\end{align*}

Therefore,
\begin{equation}\label{osamnaest}
\begin{aligned}
\sum\limits_{\substack{\alpha\in R\br{T}\\ \RE\br{\alpha}<\rho}}\RES_{s=\alpha}\br{-\frac{Z_{R}^{'}\br{s}}{Z_{R}\br{s}}\frac{x^{s+j}}{\prod\limits_{k=0}^{j}\br{s+k}}}
=O\br{x^{\rho+j}T^{-j-1+n}\log T}.
\end{aligned}
\end{equation}

It is clear that
\begin{align*}
\frac{1}{2\pi\II{}}\int\limits_{c+\II{}\tilde{T}}^{c+\II{}\infty}-\frac{Z_{R}^{'}\br{s}}{Z_{R}\br{s}}\frac{x^{s+j}}{\prod\limits_{k=0}^{j}\br{s+k}}ds
=O\br{x^{c+j}T^{-j}}
\end{align*}
and
\begin{align*}
\frac{1}{2\pi\II{}}\int\limits_{c-\II{}\infty}^{c-\II{}\tilde{T}}-\frac{Z_{R}^{'}\br{s}}{Z_{R}\br{s}}\frac{x^{s+j}}{\prod\limits_{k=0}^{j}\br{s+k}}ds
=O\br{x^{c+j}T^{-j}}.
\end{align*}

Hence,
\begin{equation}\label{devetnaest}
\begin{aligned}
\frac{1}{2\pi\II{}}\int\limits_{c-\II{}\tilde{T}}^{c+\II{}\tilde{T}}-\frac{Z_{R}^{'}\br{s}}{Z_{R}\br{s}}\frac{x^{s+j}}{\prod\limits_{k=0}^{j}\br{s+k}}ds
=\psi_{j}\br{x}+O\br{x^{c+j}T^{-j}}.
\end{aligned}
\end{equation}

Taking into account (\ref{devet})-(\ref{dvanaest}), (\ref{petnaest}), (\ref{sedamnaest})-(\ref{devetnaest}), and passing to the limit $T$ $\rightarrow$ $\infty$ in (\ref{pet}), we end up with
\begin{align*}
\psi_{j}\br{x}=&\sum\limits_{p=0}^{n-1}\br{-1}^{p+1}\sum\limits_{\br{\tau,\lambda}\in I_{p}}\sum\limits_{\rho<\alpha\leq 2\rho}\frac{x^{\alpha+j}}{\prod\limits_{k=0}^{j}\br{\alpha+k}}+\\
&\sum\limits_{p=0}^{n-1}\br{-1}^{p+1}\sum\limits_{\substack{\br{\tau,\lambda}\in I_{p}\\ \lambda=2\rho}}\sum\limits_{\RE\br{\alpha}=\rho}\frac{x^{\alpha+j}}{\prod\limits_{k=0}^{j}\br{\alpha+k}}.
\end{align*}

This completes the proof.
\end{proof}

We remark that in the case $j$ $=$ $n-1$, the $O$-terms in the relations (\ref{devet})-(\ref{dvanaest}) and (\ref{devetnaest}) all vanish during the process $T$ $\rightarrow$ $\infty$.

Notice, however, that the $O$-term in (\ref{petnaest}) is $O\br{x^{\rho+n-1}\log T}$ for $j$ $=$ $n-1$. Hence, it does not tend to zero as $T$ approaches to infinity.

This means that the method we made use of in the proof of Theorem \ref{Th1} cannot be applied in the same form to obtain the formula analogous to the formula (\ref{cetiri}) for the function $\psi_{n-1}\br{x}$.

\subsection{Explicit formulas for $\psi_{0}\br{x}$ I}
We derive the asymptotics of $\psi_{0}\br{x}$ from the asymptotics of $\psi_{j}\br{x}$, $j$ $\geq$ $n$ (see, Theorem \ref{Th1}).
\begin{theorem}\label{Th2}
Let $Y$ be as above. If $j$ $\geq$ $n$, then, for $\varepsilon$ $>$ $0$, there exists a set $E$ of finite logarithmic measure, such that
\begin{align*}
\psi_{0}\br{x}=&\sum\limits_{p=0}^{n-1}\br{-1}^{p+1}\sum\limits_{\br{\tau,\lambda}\in I_{p}}
\sum\limits_{2\rho-\rho\frac{2j+1}{2nj+1}<\alpha\leq2\rho}\frac{x^{\alpha}}{\alpha}+\\
&O\br{x^{2\rho-\rho\frac{2j+1}{2nj+1}}\br{\log x}^{\frac{n-1}{2nj+1}}\br{\log\log x}^{\frac{n-1}{2nj+1}+\varepsilon}},
\end{align*}
as $x$ $\rightarrow$ $\infty$, $x$ $\in$ $\mathbb{R}$ $\backslash$ $E$, where $\alpha$ is a singularity of the Selberg zeta function $Z_{S}\br{s+\rho-\lambda,\tau}$.
\end{theorem}
\begin{proof}
We write the second sum on the right-hand side of (\ref{cetiri}) as follows
\begin{equation}\label{dvadeset}
\begin{aligned}
&\sum\limits_{p=0}^{n-1}\br{-1}^{p+1}\sum\limits_{\substack{\br{\tau,\lambda}\in I_{p}\\ \lambda=2\rho}}
\sum\limits_{\substack{\RE\br{\alpha}=\rho\\ \abs{\IM\br{\alpha}}\leq Y}}\frac{x^{\alpha+j}}{\prod\limits_{k=0}^{j}\br{\alpha+k}}+\\
&\sum\limits_{p=0}^{n-1}\br{-1}^{p+1}\sum\limits_{\substack{\br{\tau,\lambda}\in I_{p}\\ \lambda=2\rho}}
\sum\limits_{\substack{\RE\br{\alpha}=\rho\\ Y<\abs{\IM\br{\alpha}}\leq W}}\frac{x^{\alpha+j}}{\prod\limits_{k=0}^{j}\br{\alpha+k}}+\\
&\sum\limits_{p=0}^{n-1}\br{-1}^{p+1}\sum\limits_{\substack{\br{\tau,\lambda}\in I_{p}\\ \lambda=2\rho}}
\sum\limits_{\substack{\RE\br{\alpha}=\rho\\ \abs{\IM\br{\alpha}}>W}}\frac{x^{\alpha+j}}{\prod\limits_{k=0}^{j}\br{\alpha+k}}.
\end{aligned}
\end{equation}

Following lines of \cite[p. 10]{MDPI1}, we focus our attention to the second sum in (\ref{dvadeset}).

Let $E_{p,\tau}^{i}$ denote the set
\begin{align*}
\set{x\in\sppbr{e^{i},e^{i+1}}\,:\,\abs{\sum\limits_{\substack{\RE\br{\alpha}=\rho\\ Y<\abs{\IM\br{\alpha}}\leq W}}\frac{x^{\alpha+j}}{\prod\limits_{k=0}^{j}\br{\alpha+k}}}>x^{\gamma}\br{\log x}^{\beta}\br{\log\log x}^{\beta+\varepsilon}}.
\end{align*}

It is understood that the element $\alpha$ which appears in the definition of the set $E_{p,\tau}^{i}$ is a singularity of the Selberg zeta function $Z_{S}\br{s+\rho-\lambda,\tau}$ corresponding to $p$ $\in$ $\set{0,1,...,n-1}$ and $\br{\tau,\lambda}$ $\in$ $I_{p}$, $\lambda$ $=$ $2\rho$.

Estimating the logarithmic measure $\mu^{\times}E_{p,\tau}^{i}$ $=$ $\int\limits_{E_{p,\tau}^{i}}\frac{dx}{x}$ of the set $E_{p,\tau}^{i}$ in the same way as in \cite[pp. 10-11]{MDPI1}, we obtain that $\mu^{\times}E_{p,\tau}^{i}$ is
\begin{equation}\label{dvadesetjedan}
\begin{aligned}
&O\br{\frac{e^{2\br{\rho+j-\gamma}i}}{i^{2\beta}\br{\log i}^{2\beta+2\varepsilon}}
\int\limits_{e^{i}}^{e^{i+1}}
\abs{\sum\limits_{\substack{\RE\br{\alpha}=\rho\\ Y<\abs{\IM\br{\alpha}}\leq W}}
\frac{x^{\II\IM\br{\alpha}}}{\prod\limits_{k=0}^{j}\br{\alpha+k}}}^{2}\frac{dx}{x}}\\
=&O\br{\frac{e^{2\br{\rho+j-\gamma}i}}{i^{2\beta}\br{\log i}^{2\beta+2\varepsilon}}
\int\limits_{-\frac{1}{4\pi}}^{\frac{1}{4\pi}}
\abs{\sum\limits_{\substack{\RE\br{\alpha}=\rho\\ Y<\abs{\IM\br{\alpha}}\leq W}}
\frac{e^{\II\IM\br{\alpha}\br{i+\frac{1}{2}}}}{\prod\limits_{k=0}^{j}\br{\alpha+k}}e^{2\pi\II\IM\br{\alpha}u}}^{2}du}.
\end{aligned}
\end{equation}

The Gallagher lemma \cite[p. 78, Lemma 1]{9MDPI1} implies that the integral on the right-hand side of (\ref{dvadesetjedan}) is dominated by
\begin{equation}\label{dvadesetdva}
\begin{aligned}
\int\limits_{-\infty}^{+\infty}\br{\sum\limits_{\substack{t\leq\abs{\IM\br{\alpha}}\leq t+1\\Y<\abs{\IM\br{\alpha}}\leq W}}\frac{1}{\prod\limits_{k=0}^{j}\abs{\alpha+k}}}^{2}dt.
\end{aligned}
\end{equation}

However, the integral in (\ref{dvadesetdva}) is $O\br{\frac{1}{Y^{2j+3-2n}}}$ according to Weyl's asymptotic law $N_{S,p}^{\tau,\lambda}\br{t}$ $=$ $C_{1}t^{n}$ $+$ $O\br{t^{n-1}}$.

Hence, $\mu^{\times}E_{p,\tau}^{i}$ $=$ $O\br{\frac{e^{2\br{\rho+j-\gamma}i}}{Y^{2j+3-2n}i^{2\beta}\br{\log i}^{2\beta+2\varepsilon}}}$.

For $Y$ $\sim$ $e^{\frac{\br{2\rho+2j-2\gamma}i}{2j+3-2n}}i^{\frac{1-2\beta}{2j+3-2n}}\br{\log i}^{\frac{1-2\beta}{2j+3-2n}}$, it follows that\\ $\mu^{\times}E_{p,\tau}^{i}$ $=$ $O\br{\frac{1}{i\br{\log i}^{1+2\varepsilon}}}$.

Consequently, $\mu^{\times}\bigcup\limits_{p}\bigcup\limits_{\substack{\br{\tau,\lambda}\in I_{p}\\ \lambda=2\rho}}E_{p,\tau}^{i}$ $=$ $O\br{\frac{1}{i\br{\log i}^{1+2\varepsilon}}}$.

Let $E$ denote the set $\bigcup\limits_{i}\bigcup\limits_{p}\bigcup\limits_{\substack{\br{\tau,\lambda}\in I_{p}\\ \lambda=2\rho}}E_{p,\tau}^{i}$.

The logarithmic measure of the set $E$ is finite according to the Maclaurin-Cauchy test (Cf. \cite[p. 80]{9MDPI1})
\begin{align*}
\mu^{\times}E\ll\sum\limits_{i}\frac{1}{i\br{\log i}^{1+2\varepsilon}}<\infty.
\end{align*}

Now, by the definition of $E$, it holds for any $x$ $\in$ $\mathbb{R}$ $\backslash$ $E$ that
\begin{align*}
\sum\limits_{p=0}^{n-1}\br{-1}^{p+1}\sum\limits_{\substack{\br{\tau,\lambda}\in I_{p}\\ \lambda=2\rho}}\sum\limits_{\substack{\RE\br{\alpha}=\rho\\ Y<\abs{\IM\br{\alpha}}\leq W}}\frac{x^{\alpha+j}}{\prod\limits_{k=0}^{j}\br{\alpha+k}}=O\br{x^{\gamma}\br{\log x}^{\beta}\br{\log\log x}^{\beta+\varepsilon}}.
\end{align*}

Therefore, for $x$ $\in$ $\mathbb{R}$ $\backslash$ $E$
\begin{equation}\label{dvadesettri}
\begin{aligned}
&d^{-j}\Delta_{j}^{+}\sum\limits_{p=0}^{n-1}\br{-1}^{p+1}\sum\limits_{\substack{\br{\tau,\lambda}\in I_{p}\\ \lambda=2\rho}}\sum\limits_{\substack{\RE\br{\alpha}=\rho\\ Y<\abs{\IM\br{\alpha}}\leq W}}\frac{x^{\alpha+j}}{\prod\limits_{k=0}^{j}\br{\alpha+k}}\\
=&O\br{\frac{x^{\gamma}\br{\log x}^{\beta}\br{\log\log x}^{\beta+\varepsilon}}{d^{j}}},
\end{aligned}
\end{equation}
where the operator $\Delta_{j}^{+}$ is introduced classically by
\begin{align*}
\Delta_{j}^{+} f\br{x}=\int\limits_{x}^{x+d}\int\limits_{t_{j}}^{t_{j}+d}...\int\limits_{t_{2}}^{t_{2}+d}f^{\br{j}}\br{t_{1}}dt_{1}...dt_{j}
\end{align*}
for at least $j$ times differentiable function $f$, and a constant $d$ which will be fixed later.

However, we require the bound $d$ $=$ $O\br{x}$.

It is easily seen that the application of the operator $d^{-j}\Delta_{j}^{+}$ to the first and the second sum in (\ref{dvadeset}) yields the bounds (Cf. \cite[p. 8, equation (11)]{MDPI1})
\begin{equation}\label{dvadesetcetiri}
\begin{aligned}
O\br{x^{\rho}Y^{n-1}}
\end{aligned}
\end{equation}
and
\begin{equation}\label{dvadesetpet}
\begin{aligned}
O\br{\frac{x^{\rho+j}}{d^{j}W^{j+1-n}}},
\end{aligned}
\end{equation}
respectively.

Since $\Delta x^{r}$ $=$ $d^{2n}r\br{r-1}...\br{r-\br{2n-1}}\tilde{x}^{r-2n}$ for some $\tilde{x}$ $\in$ $\SqBr{x,x+jd}$, the first sum in (\ref{cetiri}) gives us
\begin{equation}\label{dvadesetsest}
\begin{aligned}
&\sum\limits_{p=0}^{n-1}\br{-1}^{p+1}\sum\limits_{\br{\tau,\lambda}\in I_{p}}\sum\limits_{\rho<\alpha\leq2\rho}
d^{-j}\Delta_{j}^{+}\frac{x^{\alpha+j}}{\prod\limits_{k=0}^{j}\br{\alpha+k}}\\
=&\sum\limits_{p=0}^{n-1}\br{-1}^{p+1}\sum\limits_{\br{\tau,\lambda}\in I_{p}}\sum\limits_{\rho<\alpha\leq2\rho}
\frac{x^{\alpha}}{\alpha}+O\br{x^{2\rho-1}d}.
\end{aligned}
\end{equation}

Now, combining the relations (\ref{cetiri}), (\ref{dvadeset}), (\ref{dvadesettri})-(\ref{dvadesetsest}) together with the fact $\psi_{0}\br{x}$ $\leq$ $d^{-j}\Delta_{j}^{+}\psi_{j}\br{x}$, we end up with
\begin{equation}\label{dvadesetsedam}
\begin{aligned}
\psi_{0}\br{x}\leq &\sum\limits_{p=0}^{n-1}\br{-1}^{p+1}\sum\limits_{\br{\tau,\lambda}\in I_{p}}\sum\limits_{\rho<\alpha\leq2\rho}\frac{x^{\alpha}}{\alpha}+O\br{x^{2\rho-1}d}+
O\br{x^{\rho}Y^{n-1}}\\
&+O\br{\frac{x^{\gamma}\br{\log x}^{\beta}\br{\log\log x}^{\beta+\varepsilon}}{d^{j}}}+O\br{\frac{x^{\rho+j}}{d^{j}W^{j+1-n}}}
\end{aligned}
\end{equation}
for $x$ $\in$ $\mathbb{R}$ $\backslash$ $E$.

Following lines of \cite[p. 12]{MDPI1}, we compare the error terms $O\br{x^{2\rho-1}d}$, $O\br{x^{\rho}Y^{n-1}}$ and $O\br{\frac{x^{\gamma}\br{\log x}^{\beta}\br{\log\log x}^{\beta+\varepsilon}}{d^{j}}}$ on the right hand side of (\ref{dvadesetsedam}).

Hence, we require that $d$ $=$ $x^{1-\rho}Y^{n-1}$, while the equality $O\br{x^{2\rho-1}d}$ $=$\\ $O\br{\frac{x^{\gamma}\br{\log x}^{\beta}\br{\log\log x}^{\beta+\varepsilon}}{d^{j}}}$ is comfortably satisfied for $x^{2\rho-1}d$ $=$ $\frac{x^{\gamma}\br{\log x}^{\beta}\br{\log\log x}^{\beta}}{d^{j}}$, that is, for $d$ $=$ $x^{\frac{\gamma-2\rho+1}{j+1}}\br{\log x}^{\frac{\beta}{j+1}}\br{\log\log x}^{\frac{\beta}{j+1}}$.

Comparing the corresponding exponents of $x$ and $\log x$ (with respect to our selection $Y$ $\sim$ $x^{\frac{2\rho+2j-2\gamma}{2j+3-2n}}\br{\log x}^{\frac{1-2\beta}{2j+3-2n}}\br{\log\log x}^{\frac{1-2\beta}{2j+3-2n}}$), we obtain the following system of linear equations:
\begin{align*}
\frac{\gamma-2\rho+1}{j+1}=&1-\rho+\br{n-1}\frac{2\rho+2j-2\gamma}{2j+3-2n},\\
\frac{\beta}{j+1}=&\br{n-1}\frac{1-2\beta}{2j+3-2n}.
\end{align*}

The solution is given by:
\begin{align*}
\gamma=&\frac{2\br{n-\rho}j^{2}+\br{4n-3}\rho j+\rho+j}{2nj+1},\\
\beta=&\frac{\br{n-1}\br{j+1}}{2nj+1}.
\end{align*}

Hence,
\begin{align*}
d=&x^{1-\rho\frac{2j+1}{2nj+1}}\br{\log x}^{\frac{n-1}{2nj+1}}\br{\log\log x}^{\frac{n-1}{2nj+1}},\\
Y\sim&x^{\frac{2\rho j}{2nj+1}}\br{\log x}^{\frac{1}{2nj+1}}\br{\log\log x}^{\frac{1}{2nj+1}}.
\end{align*}

Substituting the obtained $d$ and $Y$ into (\ref{dvadesetsedam}), we arrive at
\begin{align*}
\psi_{0}\br{x}\leq&\sum\limits_{p=0}^{n-1}\br{-1}^{p+1}\sum\limits_{\br{\tau,\lambda}\in I_{p}}\sum\limits_{\rho<\alpha\leq 2\rho}\frac{x^{\alpha}}{\alpha}+\\
&O\br{x^{2\rho-\rho\frac{2j+1}{2nj+1}}\br{\log x}^{\frac{n-1}{2nj+1}}\br{\log\log x}^{\frac{n-1}{2nj+1}+\varepsilon}}\\
=&\sum\limits_{p=0}^{n-1}\br{-1}^{p+1}\sum\limits_{\br{\tau,\lambda}\in I_{p}}\sum\limits_{2\rho-\rho\frac{2j+1}{2nj+1}<\alpha\leq2\rho}\frac{x^{\alpha}}{\alpha}+\\
&O\br{x^{2\rho-\rho\frac{2j+1}{2nj+1}}\br{\log x}^{\frac{n-1}{2nj+1}}\br{\log\log x}^{\frac{n-1}{2nj+1}+\varepsilon}}
\end{align*}
as $x$ $\rightarrow$ $\infty$, $x$ $\in$ $\mathbb{R}$ $\backslash$ $E$.

In a similar way,
\begin{align*}
\psi_{0}\br{x}\geq&\sum\limits_{p=0}^{n-1}\br{-1}^{p+1}\sum\limits_{\br{\tau,\lambda}\in I_{p}}\sum\limits_{2\rho-\rho\frac{2j+1}{2nj+1}<\alpha\leq2\rho}\frac{x^{\alpha}}{\alpha}+\\
&O\br{x^{2\rho-\rho\frac{2j+1}{2nj+1}}\br{\log x}^{\frac{n-1}{2nj+1}}\br{\log\log x}^{\frac{n-1}{2nj+1}+\varepsilon}}
\end{align*}
as $x$ $\rightarrow$ $\infty$, $x$ $\in$ $\mathbb{R}$ $\backslash$ $E$.

This completes the proof.
\end{proof}

\subsection{Gallagherian prime geodesic theorem I}
An immediate consequence of Theorem \ref{Th2} is the following theorem (Cf. \cite[p. 102]{Park}).
\begin{theorem}\label{Th3}\,(Gallagherian Prime Geodesic Theorem)\,
Let $Y$ be as above. If $j$ $\geq$ $n$, then, for $\varepsilon$ $>$ $0$, there is a set $E$ of finite logarithmic measure, so that
\begin{equation}\label{dvadesetosam}
\begin{aligned}
\pi_{\Gamma}\br{x}=&\sum\limits_{p=0}^{n-1}\br{-1}^{p+1}\sum\limits_{\br{\tau,\lambda}\in I_{p}}
\sum\limits_{2\rho-\rho\frac{2j+1}{2nj+1}<\alpha\leq2\rho}\li\br{x^{\alpha}}+\\
&O\br{x^{2\rho-\rho\frac{2j+1}{2nj+1}}\br{\log x}^{\frac{n-1}{2nj+1}-1}\br{\log\log x}^{\frac{n-1}{2nj+1}+\varepsilon}},
\end{aligned}
\end{equation}
as $x$ $\rightarrow$ $\infty$, $x$ $\in$ $\mathbb{R}$ $\backslash$ $E$, where $\alpha$ is a singularity of the Selberg zeta function $Z_{S}\br{s+\rho-\lambda,\tau}$.
\end{theorem}
\vspace{2mm}

Let's have a closer look at the exponent $2\rho$ $-$ $\rho\frac{2j+1}{2nj+1}$ of $x$ in the $O$-term of Theorem \ref{Th3}.

It is understood that we assume that $n$ $\geq$ $2$.

Hence, the inequality $2\rho$ $-$ $\rho\frac{2j+1}{2nj+1}$ $<$ $2\rho$ $-$ $\rho\frac{2\br{j+1}+1}{2n\br{j+1}+1}$ is valid being equivalent to the inequality $2n$ $+$ $1$ $>$ $3$ which is obviously correct.

Thus, the sequence $\set{2\rho-\rho\frac{2j+1}{2nj+1}}_{j}$ is strictly increasing, with the limit $2\rho$ $-$ $\frac{\rho}{n}$. This means that the optimal size of the error term in (\ref{dvadesetosam}) is\\ $O\br{x^{2\rho-\rho\frac{2n+1}{2n^{2}+1}}\br{\log x}^{\frac{n-1}{2n^{2}+1}-1}\br{\log\log x}^{\frac{n-1}{2n^{2}+1}+\varepsilon}}$, and is achieved for $j$ $=$ $n$. Note that the limit $2\rho$ $-$ $\frac{\rho}{n}$ is nothing else but the optimal exponent of $x$ in DeGeorge's prime geodesic theorem \cite[p. 2, relation (1)]{MDPI1}. It is also clear that the $O$-term in (\ref{dvadesetosam}) improves the $O$-term in Theorem \ref{Th3} \cite[p. 9]{MDPI1} for $n$ $\leq$ $j$ $<$ $2n$.

Since the sequence $\set{2\rho-\rho\frac{2j+1}{2nj+1}}_{j}$ is strictly increasing, and its terms for $j$ $\geq$ $n$ are obtained by the application of the explicit formulas for the counting functions $\psi_{j}\br{x}$, $j$ $\geq$ $n$ (Theorem \ref{Th1}), the question that naturally arises is the possibility of further reduction of the error term in (\ref{dvadesetosam}) via the function $\psi_{n-1}\br{x}$. As we already noted in the remark at the end of Subsection 3.1, the method applied in the proof of Theorem \ref{Th1} cannot be copied to give us the desired form (\ref{cetiri}) for $\psi_{n-1}\br{x}$.

However, recall that Koyama \cite[p. 79, relation (5)]{9MDPI1}, applied one of Hejhal's explicit formulas for $\psi_{1}\br{x}$ (see, Theorem 6.16 \cite[p. 110]{Hejhal1} and proof of Theorem 3.4 \cite[p. 474]{Hejhal2}), to obtain the Gallagherian prime geodesic theorem for compact Riemann surfaces and generic hyperbolic surfaces of finite volume ($n$ $=$ $2$, $\rho$ $=$ $\frac{1}{2}\br{n-1}$ $=$ $\frac{1}{2}$). For the same reason, we consider $\psi_{n-1}\br{x}$ in the next section.

\subsection{Explicit formulas for $\psi_{n-1}\br{x}$}
Following Hejhal's approach in the compact Riemann surface case \cite[pp. 104-110]{Hejhal1}, we consider the counting function $\psi_{n-1}\br{x}$.

We obtain the following conditional result.
\begin{theorem}\label{Th4}
Let $Y$ be as above, and $\varepsilon_{1}$ $>$ $0$, $\delta$ $>$ $0$. Suppose $T$ $\gg$ $0$ and $A$ $\gg$ $2\rho$ are chosen so that $\II{}T$ is not a singularity of $Z_{S}\br{s,\tau}$, $\tau$ $\in$ $\mathcal{T}$, and $-A$ is not a pole of $-\frac{Z_{R}^{'}\br{s}}{Z_{R}\br{s}}$. If the integral $\int\limits_{-A+\II{}T}^{\rho+\II{}T}-\frac{Z_{R}^{'}\br{s}}{Z_{R}\br{s}}\frac{x^{s+n-1}}{\prod\limits_{k=0}^{n-1}\br{s+k}}ds$ is $O\br{\frac{x^{2\rho+\varepsilon_{1}+n-1}}{\varepsilon_{1}T^{1-\delta}}}$, then
\begin{equation}\label{dvadesetdevet}
\begin{aligned}
\psi_{n-1}\br{x}=&\sum\limits_{p=0}^{n-1}\alpha_{n-1-k}x^{n-1-k}\log x+\sum\limits_{k=0}^{n-1}\beta_{n-1-k}x^{n-1-k}+\\
&\sum\limits_{p=0}^{n-1}\br{-1}^{p+1}\sum\limits_{\br{\tau,\lambda}\in I_{p}}\sum\limits_{\alpha\in S_{p,\tau,\lambda}^{\mathbb{R}}}\frac{x^{\alpha+n-1}}{\prod\limits_{k=0}^{n-1}\br{\alpha+k}}+\\
&\sum\limits_{p=0}^{n-1}\br{-1}^{p+1}\sum\limits_{\br{\tau,\lambda}\in I_{p}}\sum\limits_{\substack{\alpha\in S_{p,\tau,\lambda}^{-\rho+\lambda}\\ \abs{\IM\br{\alpha}}\leq T}}\frac{x^{\alpha+n-1}}{\prod\limits_{k=0}^{n-1}\br{\alpha+k}}+O\br{\frac{x^{2\rho+\varepsilon_{1}+n-1}}{\varepsilon_{1}T^{1-\delta}}},
\end{aligned}
\end{equation}
where $S_{p,\tau,\lambda}^{\mathbb{R}}$ is the set of real singularities of $Z_{S}\br{s+\rho-\lambda,\tau}$ not containing the integers $0$, $-1$, ..., $-\br{n-1}$, $S_{p,\tau,\lambda}^{-\rho+\lambda}$ is the set of non-real singularities of $Z_{S}\br{s+\rho-\lambda,\tau}$, and $\alpha_{k}$, $\beta_{k}$, $k$ $\in$ $\set{0,1,...,n-1}$ are some explicitly computable constants.
\end{theorem}
\begin{proof}
By the Cauchy residue theorem applied to the rectangle $R\br{A,T}$ given by vertices: $2\rho$ $+$ $\varepsilon_{1}$ $-$ $\II{}T$, $2\rho$ $+$ $\varepsilon_{1}$ $+$ $\II{}T$, $-A$ $+$ $\II{}T$, $-A$ $-$ $\II{}T$, we immediately have that
\begin{equation}\label{trideset}
\begin{aligned}
&\frac{1}{2\pi\II{}}\int\limits_{2\rho+\varepsilon_{1}-\II{}T}^{2\rho+\varepsilon_{1}+\II{}T}-\frac{Z_{R}^{'}\br{s}}{Z_{R}\br{s}}
\frac{x^{s+n-1}}{\prod\limits_{k=0}^{n-1}\br{s+k}}ds\\
=&\frac{1}{2\pi\II{}}\int\limits_{\rho+\varepsilon_{1}+\II{}T}^{2\rho+\varepsilon_{1}+\II{}T}+
\frac{1}{2\pi\II{}}\int\limits_{\rho+\II{}T}^{\rho+\varepsilon_{1}+\II{}T}+
\frac{1}{2\pi\II{}}\int\limits_{-A+\II{}T}^{\rho+\II{}T}+
\frac{1}{2\pi\II{}}\int\limits_{-A-\II{}T}^{-A+\II{}T}\\
&-\frac{1}{2\pi\II{}}\int\limits_{-A-\II{}T}^{\rho-\II{}T}-\frac{1}{2\pi\II{}}\int\limits_{\rho-\II{}T}^{\rho+\varepsilon_{1}-\II{}T}
-\frac{1}{2\pi\II{}}\int\limits_{\rho+\varepsilon_{1}-\II{}T}^{2\rho+\varepsilon_{1}-\II{}T}+
\end{aligned}
\end{equation}
\begin{align*}
&\sum\limits_{\alpha\in R\br{A,T}}\RES_{s=\alpha}\br{-\frac{Z_{R}^{'}\br{s}}{Z_{R}\br{s}}\frac{x^{s+n-1}}{\prod\limits_{k=0}^{n-1}\br{s+k}}}.
\end{align*}

As in the derivation of (\ref{devetnaest}), we obtain
\begin{equation}\label{tridesetjedan}
\begin{aligned}
\frac{1}{2\pi\II{}}\int\limits_{2\rho+\varepsilon_{1}-\II{}T}^{2\rho+\varepsilon_{1}+\II{}T}-\frac{Z_{R}^{'}\br{s}}{Z_{R}\br{s}}
\frac{x^{s+n-1}}{\prod\limits_{k=0}^{n-1}\br{s+k}}ds=\psi_{n-1}\br{x}+O\br{\frac{x^{2\rho+\varepsilon_{1}+n-1}}{T^{n-1}}}.
\end{aligned}
\end{equation}

By (\ref{osam}),
\begin{equation}\label{tridesetdva}
\begin{aligned}
\frac{1}{2\pi\II{}}\int\limits_{\rho+\varepsilon_{1}+\II{}T}^{2\rho+\varepsilon_{1}+\II{}T}-\frac{Z_{R}^{'}\br{s}}{Z_{R}\br{s}}
\frac{x^{s+n-1}}{\prod\limits_{k=0}^{n-1}\br{s+k}}ds=O\br{\frac{x^{2\rho+\varepsilon_{1}+n-1}}{\varepsilon_{1}T^{1-\delta}}}
\end{aligned}
\end{equation}
and
\begin{equation}\label{tridesettri}
\begin{aligned}
-\frac{1}{2\pi\II{}}\int\limits_{\rho+\varepsilon_{1}-\II{}T}^{2\rho+\varepsilon_{1}-\II{}T}-\frac{Z_{R}^{'}\br{s}}{Z_{R}\br{s}}
\frac{x^{s+n-1}}{\prod\limits_{k=0}^{n-1}\br{s+k}}ds=O\br{\frac{x^{2\rho+\varepsilon_{1}+n-1}}{\varepsilon_{1}T^{1-\delta}}}.
\end{aligned}
\end{equation}

Reasoning in the same way as in the derivation of (\ref{devet}) and (\ref{deset}), we have
\begin{equation}\label{tridesetcetiri}
\begin{aligned}
\frac{1}{2\pi\II{}}\int\limits_{\rho+\II{}T}^{\rho+\varepsilon_{1}+\II{}T}-\frac{Z_{R}^{'}\br{s}}{Z_{R}\br{s}}
\frac{x^{s+n-1}}{\prod\limits_{k=0}^{n-1}\br{s+k}}ds=O\br{\frac{x^{\rho+\varepsilon_{1}+n-1}}{T^{1-\delta}}}
\end{aligned}
\end{equation}
and
\begin{equation}\label{tridesetpet}
\begin{aligned}
-\frac{1}{2\pi\II{}}\int\limits_{\rho-\II{}T}^{\rho+\varepsilon_{1}-\II{}T}-\frac{Z_{R}^{'}\br{s}}{Z_{R}\br{s}}
\frac{x^{s+n-1}}{\prod\limits_{k=0}^{n-1}\br{s+k}}ds=O\br{\frac{x^{\rho+\varepsilon_{1}+n-1}}{T^{1-\delta}}}.
\end{aligned}
\end{equation}

Note that the $O$-terms in (\ref{tridesetjedan})-(\ref{tridesetpet}) are all dominated by $O\br{\frac{x^{2\rho+\varepsilon_{1}+n-1}}{\varepsilon_{1}T^{1-\delta}}}$.

We calculate the residues in (\ref{trideset}) in the same way as in \cite[pp. 6-7]{MDPI1}. We obtain,
\begin{equation}\label{tridesetsest}
\begin{aligned}
&\sum\limits_{\alpha\in R\br{A,T}}\RES_{s=\alpha}\br{-\frac{Z_{R}^{'}\br{s}}{Z_{R}\br{s}}\frac{x^{s+n-1}}{\prod\limits_{k=0}^{n-1}\br{s+k}}}\\
=&\sum\limits_{k=0}^{n-1}\alpha_{n-1-k}x^{n-1-k}\log x+\sum\limits_{k=0}^{n-1}\beta_{n-1-k}x^{n-1-k}+\\
&\sum\limits_{p=0}^{n-1}\br{-1}^{p+1}\sum\limits_{\br{\tau,\lambda}\in I_{p}}\sum\limits_{\substack{\alpha\in S_{p,\tau,\lambda}^{\mathbb{R}}\\ \alpha>-A}}
\frac{x^{\alpha+n-1}}{\prod\limits_{k=0}^{n-1}\br{\alpha+k}}+\\
&\sum\limits_{p=0}^{n-1}\br{-1}^{p+1}\sum\limits_{\br{\tau,\lambda}\in I_{p}}\sum\limits_{\substack{\alpha\in S_{p,\tau,\lambda}^{-\rho+\lambda}\\ \abs{\IM\br{\alpha}}\leq T}}
\frac{x^{\alpha+n-1}}{\prod\limits_{k=0}^{n-1}\br{\alpha+k}}.
\end{aligned}
\end{equation}

Combining the relations (\ref{trideset})-(\ref{tridesetsest}) together with the assumption\\ $\int\limits_{-A+\II{}T}^{\rho+\II{}T}-\frac{Z_{R}^{'}\br{s}}{Z_{R}\br{s}}\frac{x^{s+n-1}}{\prod\limits_{k=0}^{n-1}\br{s+k}}ds$ $=$ $O\br{\frac{x^{2\rho+\varepsilon_{1}+n-1}}{\varepsilon_{1}T^{1-\delta}}}$, we arrive at

\begin{align*}
\psi_{n-1}\br{x}=&\sum\limits_{k=0}^{n-1}\alpha_{n-1-k}x^{n-1-k}\log x+\sum\limits_{k=0}^{n-1}\beta_{n-1-k}x^{n-1-k}+\\
&\sum\limits_{p=0}^{n-1}\br{-1}^{p+1}\sum\limits_{\br{\tau,\lambda}\in I_{p}}\sum\limits_{\substack{\alpha\in S_{p,\tau,\lambda}^{\mathbb{R}}\\ \alpha>-A}}
\frac{x^{\alpha+n-1}}{\prod\limits_{k=0}^{n-1}\br{\alpha+k}}+
\end{align*}
\begin{align*}
&\sum\limits_{p=0}^{n-1}\br{-1}^{p+1}\sum\limits_{\br{\tau,\lambda}\in I_{p}}\sum\limits_{\substack{\alpha\in S_{p,\tau,\lambda}^{-\rho+\lambda}\\ \abs{\IM\br{\alpha}}\leq T}}
\frac{x^{\alpha+n-1}}{\prod\limits_{k=0}^{n-1}\br{\alpha+k}}+\\
&O\br{\frac{x^{2\rho+\varepsilon_{1}+n-1}}{\varepsilon_{1}T^{1-\delta}}}+
\frac{1}{2\pi\II{}}\int\limits_{-A-\II{}T}^{-A+\II{}T}-\frac{Z_{R}^{'}\br{s}}{Z_{R}\br{s}}
\frac{x^{s+n-1}}{\prod\limits_{k=0}^{n-1}\br{s+k}}ds.
\end{align*}

Passing to the limit $A$ $\rightarrow$ $+\infty$, (Cf. also, \cite[p. 85, Prop. 5.7]{Hejhal1} and \cite[p. 244]{Randol}), we obtain the assertion.

This completes the proof.
\end{proof}

\subsection{Explicit formulas for $\psi_{0}\br{x}$ II}
Now, we derive the asymptotics of $\psi_{0}\br{x}$ from the asymptotics of $\psi_{n-1}\br{x}$ (see, Theorem \ref{Th4}).
\begin{theorem}\label{Th5}
Suppose that the assumptions of Theorem \ref{Th4} hold. For $\varepsilon$ $>$ $0$, there is a set $E$ of finite logarithmic measure, such that
\begin{align*}
\psi_{0}\br{x}=&\sum\limits_{p=0}^{n-1}\br{-1}^{p+1}\sum\limits_{\br{\tau,\lambda}\in I_{p}}
\sum\limits_{\substack{\alpha\in S_{p,\tau,\lambda}^{\mathbb{R}}\\2\rho-\rho\frac{2\br{n-1}+1}{2n\br{n-1}+1}<\alpha\leq2\rho}}\frac{x^{\alpha}}{\alpha}+\\
&O\br{x^{2\rho-\rho\frac{2\br{n-1}+1}{2n\br{n-1}+1}}\br{\log x}^{\frac{n-1}{2n\br{n-1}+1}}\br{\log\log x}^{\frac{n-1}{2n\br{n-1}+1}+\varepsilon}},
\end{align*}
as $x$ $\rightarrow$ $\infty$, $x$ $\in$ $\mathbb{R}$ $\backslash$ $E$.
\end{theorem}
\begin{proof}
We write the third sum on the right-hand side of (\ref{dvadesetdevet}) as follows
\begin{equation}\label{tridesetsedam}
\begin{aligned}
&\sum\limits_{p=0}^{n-1}\br{-1}^{p+1}\sum\limits_{\br{\tau,\lambda}\in I_{p}}\sum\limits_{\substack{\alpha\in S_{p,\tau,\lambda}^{-\rho+\lambda}\\ \abs{\IM\br{\alpha}}\leq Y}}
\frac{x^{\alpha+n+1}}{\prod\limits_{k=0}^{n-1}\br{\alpha+k}}+\\
&\sum\limits_{p=0}^{n-1}\br{-1}^{p+1}\sum\limits_{\br{\tau,\lambda}\in I_{p}}\sum\limits_{\substack{\alpha\in S_{p,\tau,\lambda}^{-\rho+\lambda}\\ Y<\abs{\IM\br{\alpha}}\leq T}}\frac{x^{\alpha+n+1}}{\prod\limits_{k=0}^{n-1}\br{\alpha+k}}.
\end{aligned}
\end{equation}

Reasoning in the same way as in the proof of Theorem \ref{Th3} in \cite[pp. 10-11]{MDPI1}, we obtain that $\mu^{\times}E_{p,\tau,\lambda}^{i}$ $=$ $\int\limits_{E_{p,\tau,\lambda}^{i}}\frac{dx}{x}$ $=$ $O\br{\frac{e^{2\br{\rho+n-1-\gamma}i}}{Yi^{2\beta}\br{\log i}^{2\beta+2\varepsilon}}}$, where the sets $E_{p,\tau,\lambda}^{i}$ have the same form as before, with $j$ replaced by $i$, $2n$ replaced with $n-1$, $W$ $=$ $T$, and $\alpha$ in the expression $x^{\alpha}\br{\log x}^{\beta}\br{\log\log x}^{\beta+\varepsilon}$ replaced with $\gamma$.

The relations (15) and (16) in \cite{MDPI1} are adjusted with respect to adapted notation.

In particular,
\begin{align*}
\int\limits_{-\infty}^{+\infty}\br{\sum\limits_{\substack{t\leq\abs{\IM\br{\alpha}}\leq t+1\\Y<\abs{\IM\br{\alpha}}\leq T}}\frac{1}{\prod\limits_{k=0}^{n-1}\abs{\alpha+k}}}^{2}dt=O\br{\frac{1}{Y}}.
\end{align*}

Put $Y$ $\sim$ $e^{2\br{\rho+n-1-\gamma}i}i^{1-2\beta}\br{\log i}^{1-2\beta}$.

It follows that $\mu^{\times}E$ $<$ $\infty$, where $E$ $=$ $\bigcup\limits_{i}\bigcup\limits_{p}\bigcup\limits_{\br{\tau,\lambda}\in I_{p}}E_{p,\tau,\lambda}^{i}$.

We obtain:
\begin{equation}\label{tridesetosam}
\begin{aligned}
&d^{-\br{n-1}}\Delta_{n-1}^{+}\sum\limits_{p=0}^{n-1}\br{-1}^{p+1}\sum\limits_{\br{\tau,\lambda}\in I_{p}}\sum\limits_{\substack{\alpha\in S_{p,\tau,\lambda}^{-\rho+\lambda}\\ Y<\abs{\IM\br{\alpha}}\leq T}}\frac{x^{\alpha+n-1}}{\prod\limits_{k=0}^{n-1}\br{\alpha+k}}\\
=&O\br{\frac{x^{\gamma}\br{\log x}^{\beta}\br{\log\log x}^{\beta+\varepsilon}}{d^{n-1}}}
\end{aligned}
\end{equation}
for $x$ $\in$ $\mathbb{R}$ $\backslash$ $E$,
\begin{equation}\label{tridesetdevet}
\begin{aligned}
d^{-\br{n-1}}\Delta_{n-1}^{+}\sum\limits_{p=0}^{n-1}\br{-1}^{p+1}\sum\limits_{\br{\tau,\lambda}\in I_{p}}\sum\limits_{\substack{\alpha\in S_{p,\tau,\lambda}^{-\rho+\lambda}\\ \abs{\IM\br{\alpha}}\leq Y}}\frac{x^{\alpha+n-1}}{\prod\limits_{k=0}^{n-1}\br{\alpha+k}}
=O\br{x^{\rho}Y^{n-1}},
\end{aligned}
\end{equation}
\begin{equation}\label{cetrdeset}
\begin{aligned}
&d^{-\br{n-1}}\Delta_{n-1}^{+}\br{\sum\limits_{k=0}^{n-1}\alpha_{n-1-k}x^{n-1-k}\log x+\sum\limits_{k=0}^{n-1}\beta_{n-1-k}x^{n-1-k}}\\
=&O\br{\log x},
\end{aligned}
\end{equation}
\begin{equation}\label{cetrdesetjedan}
\begin{aligned}
&d^{-\br{n-1}}\Delta_{n-1}^{+}\sum\limits_{p=0}^{n-1}\br{-1}^{p+1}\sum\limits_{\br{\tau,\lambda}\in I_{p}}\sum\limits_{\alpha\in S_{p,\tau,\lambda}^{\mathbb{R}}}
\frac{x^{\alpha+n-1}}{\prod\limits_{k=0}^{n-1}\br{\alpha+k}}\\
=&\sum\limits_{p=0}^{n-1}\br{-1}^{p+1}\sum\limits_{\br{\tau,\lambda}\in I_{p}}\sum\limits_{\substack{\alpha\in S_{p,\tau,\lambda}^{\mathbb{R}}\\ 0<\alpha\leq2\rho}}
\frac{x^{\alpha}}{\alpha}+O\br{x^{2\rho-1}d}.
\end{aligned}
\end{equation}

Combining the relations (\ref{dvadesetdevet}), (\ref{tridesetsedam})-(\ref{cetrdesetjedan}) for $x$ $\in$ $\mathbb{R}$ $\backslash$ $E$ together with the fact that\\ $\psi_{0}\br{x}$ $\leq$ $d^{-\br{n-1}}\Delta_{n-1}^{+}\psi_{n-1}\br{x}$, we end up with
\begin{equation}\label{cetrdesetdva}
\begin{aligned}
\psi_{0}\br{x}\leq&\sum\limits_{p=0}^{n-1}\br{-1}^{p+1}\sum\limits_{\br{\tau,\lambda}\in I_{p}}\sum\limits_{\substack{\alpha\in S_{p,\tau,\lambda}^{\mathbb{R}}\\ 0<\alpha\leq2\rho}}
\frac{x^{\alpha}}{\alpha}+O\br{x^{2\rho-1}d}+\\
&O\br{x^{\rho}Y^{n-1}}+O\br{\frac{x^{\gamma}\br{\log x}^{\beta}\br{\log\log x}^{\beta+\varepsilon}}{d^{n-1}}}+\\
&O\br{\frac{x^{2\rho+\varepsilon_{1}+n-1}}{\varepsilon_{1}d^{n-1}T^{1-\delta}}}.
\end{aligned}
\end{equation}

We compare the error terms $O\br{x^{2\rho-1}d}$, $O\br{x^{\rho}Y^{n-1}}$ and\\ $O\br{\frac{x^{\gamma}\br{\log x}^{\beta}\br{\log\log x}^{\beta+\varepsilon}}{d^{n-1}}}$ in the same way as earlier.

Note that it is enough to replace $j$ with $n-1$ in the corresponding relations in the proof of Theorem \ref{Th2}.

We obtain the system:
\begin{align*}
\frac{\gamma-2\rho+1}{n}=&1-\rho+2\br{n-1}\br{\rho+n-1-\gamma},
\end{align*}
\begin{align*}
\frac{\beta}{n}=&\br{n-1}\br{1-2\beta},
\end{align*}
its solution:
\begin{align*}
\gamma=&\frac{2\br{n-\rho}\br{n-1}^{2}+\br{4n-3}\br{n-1}\rho+\rho+\br{n-1}}{2n\br{n-1}+1},\\
\beta=&\frac{\br{n-1}n}{2n\br{n-1}+1},
\end{align*}
and:
\begin{align*}
d=&x^{1-\rho\frac{2\br{n-1}+1}{2n\br{n-1}+1}}\br{\log x}^{\frac{n-1}{2n\br{n-1}+1}}\br{\log\log x}^{\frac{n-1}{2n\br{n-1}+1}},\\
Y\sim&x^{\frac{2\rho\br{n-1}}{2n\br{n-1}+1}}\br{\log x}^{\frac{1}{2n\br{n-1}+1}}\br{\log\log x}^{\frac{1}{2n\br{n-1}+1}}.
\end{align*}

Substituting the obtained $d$ and $Y$ into (\ref{cetrdesetdva}), we finally have that for $x$ $\in$ $\mathbb{R}$ $\backslash$ $E$
\begin{align*}
\psi_{0}\br{x}\leq&\sum\limits_{p=0}^{n-1}\br{-1}^{p+1}\sum\limits_{\br{\tau,\lambda}\in I_{p}}\sum\limits_{\substack{\alpha\in S_{p,\tau,\lambda}^{\mathbb{R}}\\ 2\rho-\rho\frac{2\br{n-1}+1}{2n\br{n-1}+1}<\alpha\leq 2\rho}}\frac{x^{\alpha}}{\alpha}+\\
&O\br{x^{2\rho-\rho\frac{2\br{n-1}+1}{2n\br{n-1}+1}}\br{\log x}^{\frac{n-1}{2n\br{n-1}+1}}\br{\log\log x}^{\frac{n-1}{2n\br{n-1}+1}+\varepsilon}}.
\end{align*}

Since the opposite inequality holds also true, the assertion follows.

This completes the proof.
\end{proof}

Now, we are in position to improve the result given by Theorem \ref{Th3}.

\subsection{Gallagherian prime geodesic theorem II}
As in the case of Theorem \ref{Th3}, the direct consequence of Theorem \ref{Th5} is the following theorem.
\begin{theorem}\label{Th6}\,(Gallagherian Prime Geodesic Theorem)\,
Suppose that the assumptions of Theorem \ref{Th4} hold. For $\varepsilon$ $>$ $0$, there is a set $E$ of finite logarithmic measure, such that
\begin{align*}
\pi_{\Gamma}\br{x}=&\sum\limits_{p=0}^{n-1}\br{-1}^{p+1}\sum\limits_{\br{\tau,\lambda}\in I_{p}}
\sum\limits_{\substack{\alpha\in S_{p,\tau,\lambda}^{\mathbb{R}}\\2\rho-\rho\frac{2\br{n-1}+1}{2n\br{n-1}+1}<\alpha\leq2\rho}}\li\br{x^{\alpha}}+\\
&O\br{x^{2\rho-\rho\frac{2\br{n-1}+1}{2n\br{n-1}+1}}\br{\log x}^{\frac{n-1}{2n\br{n-1}+1}-1}\br{\log\log x}^{\frac{n-1}{2n\br{n-1}+1}+\varepsilon}},
\end{align*}
as $x$ $\rightarrow$ $\infty$, $x$ $\in$ $\mathbb{R}$ $\backslash$ $E$.
\end{theorem}

\subsection{Prime geodesic theorem}
In \cite{MDPI1}(see, Theorem 2), we derived the prime geodesic theorem by applying the counting function $\psi_{2n}\br{x}$.

In this section, we prove that taking the function $\psi_{j}\br{x}$ with $j$ $\geq$ $n$ (see, Theorem \ref{Th1}), does not yield a better result.

Indeed (Cf. \cite[pp. 370-371]{Koreja}), by (\ref{cetiri}),
\begin{align*}
\psi_{0}\br{x}\leq& d^{-j}\Delta_{j}^{+}\psi_{j}\br{x}\\
=&\sum\limits_{p=0}^{n-1}\br{-1}^{p+1}\sum\limits_{\br{\tau,\lambda}\in I_{p}}
\sum\limits_{\rho<\alpha\leq 2\rho}d^{-j}\Delta_{j}^{+}\frac{x^{\alpha+j}}{\prod\limits_{k=0}^{j}\br{\alpha+k}}+\\
&\sum\limits_{p=0}^{n-1}\br{-1}^{p+1}\sum\limits_{\substack{\br{\tau,\lambda}\in I_{p}\\ \lambda=2\rho}}
\sum\limits_{\RE\br{\alpha}=\rho}d^{-j}\Delta_{j}^{+}\frac{x^{\alpha+j}}{\prod\limits_{k=0}^{j}\br{\alpha+k}}\\
=&\sum\limits_{p=0}^{n-1}\br{-1}^{p+1}\sum\limits_{\br{\tau,\lambda}\in I_{p}}
\sum\limits_{\rho<\alpha\leq 2\rho}\frac{x^{\alpha}}{\alpha}+O\br{x^{2\rho-1}d}+\\
&\sum\limits_{p=0}^{n-1}\br{-1}^{p+1}\sum\limits_{\substack{\br{\tau,\lambda}\in I_{p}\\ \lambda=2\rho}}O\br{x^{\rho}\int\limits_{\rho}^{K}t^{-1}dN_{S,p}^{\tau,\lambda}\br{t}}+\\
&\sum\limits_{p=0}^{n-1}\br{-1}^{p+1}\sum\limits_{\substack{\br{\tau,\lambda}\in I_{p}\\ \lambda=2\rho}}O\br{d^{-j}x^{\rho+j}\int\limits_{K}^{+\infty}t^{-j-1}dN_{S,p}^{\tau,\lambda}\br{t}}\\
=&\sum\limits_{p=0}^{n-1}\br{-1}^{p+1}\sum\limits_{\br{\tau,\lambda}\in I_{p}}
\sum\limits_{\rho<\alpha\leq 2\rho}\frac{x^{\alpha}}{\alpha}+O\br{x^{2\rho-1}d}+O\br{x^{\rho}K^{n-1}}+\\
&O\br{d^{-j}x^{\rho+j}K^{n-j-1}}.
\end{align*}

Substituting $d$ $=$ $x^{1-\frac{\rho}{n}}$, $K$ $=$ $x^{\frac{\rho}{n}}$, we get
\begin{align*}
\psi_{0}\br{x}\leq\sum\limits_{p=0}^{n-1}\br{-1}^{p+1}\sum\limits_{\br{\tau,\lambda}\in I_{p}}
\sum\limits_{2\rho-\frac{\rho}{n}<\alpha\leq 2\rho}\frac{x^{\alpha}}{\alpha}+O\br{x^{2\rho-\frac{\rho}{n}}},
\end{align*}
that is,
\begin{equation}\label{cetrdesettri}
\begin{aligned}
\pi_{\Gamma}\br{x}=\sum\limits_{p=0}^{n-1}\br{-1}^{p+1}\sum\limits_{\br{\tau,\lambda}\in I_{p}}
\sum\limits_{2\rho-\frac{\rho}{n}<\alpha\leq 2\rho}\li\br{x^{\alpha}}+O\br{x^{2\rho-\frac{\rho}{n}}\br{\log x}^{-1}}
\end{aligned}
\end{equation}
as $x$ $\rightarrow$ $\infty$.

Note that one is in position to derive a variant of the asymptotics (\ref{cetrdesettri}) from the asymptotics of $\psi_{n-1}\br{x}$ (Cf. \cite[pp. 98-102]{Park}).

However, such approach yields only a weaker form of the prime geodesic theorem, with $O\br{x^{2\rho-\frac{\rho}{n}}\br{\log x}^{-\frac{1}{2}}}$ in place of $O\br{x^{2\rho-\frac{\rho}{n}}\br{\log x}^{-1}}$ in (\ref{cetrdesettri}).

\section{Discussion}
The error terms in (\ref{dvadesetosam}), i.e., the terms\\ $O\br{x^{2\rho-\rho\frac{2j+1}{2nj+1}}\br{\log x}^{\frac{n-1}{2nj+1}-1}\br{\log\log x}^{\frac{n-1}{2nj+1}+\varepsilon}}$, $j$ $\geq$ $n$, determine the sequence $\set{2\rho-\rho\frac{2j+1}{2nj+1}}_{j}$, whose elements are given by the exponents of $x$, and are such that $2\rho$ $-$ $\rho\frac{2j+1}{2nj+1}$ $<$ $2\rho$ $-$ $\rho\frac{2\br{j+1}+1}{2n\br{j+1}+1}$ for $j$ $\in$ $\mathbb{N}$. The sequence is strictly increasing since the last inequality is equivalent to the inequality $n$ $>$ $1$, which is obviously true. Being interested in $j$ $\geq$ $n$ in Theorem \ref{Th3}, we conclude that the optimal size of the error term in (\ref{dvadesetosam}) is $O\br{x^{2\rho-\rho\frac{2n+1}{2n^{2}+1}}\br{\log x}^{\frac{n-1}{2n^{2}+1}-1}\br{\log\log x}^{\frac{n-1}{2n^{2}+1}+\varepsilon}}$, and is achieved for $j$ $=$ $n$. So, this $O$-term, as well as each of the $O$-terms $O\br{x^{2\rho-\rho\frac{2j+1}{2nj+1}}\br{\log x}^{\frac{n-1}{2nj+1}-1}\br{\log\log x}^{\frac{n-1}{2nj+1}+\varepsilon}}$, $n$ $<$ $j$ $<$ $2n$, improve the $O$-term $O\br{x^{2\rho-\rho\frac{2\cdot\br{2n}+1}{2n\cdot\br{2n}+1}}\br{\log x}^{\frac{n-1}{2n\cdot\br{2n}+1}-1}\br{\log\log x}^{\frac{n-1}{2n\cdot\br{2n}+1}+\varepsilon}}$ derived in our previous study \cite[p. 9, Th. 3]{MDPI1}.

Note that our $O\br{x^{2\rho-\rho\frac{2n+1}{2n^{2}+1}}\br{\log x}^{\frac{n-1}{2n^{2}+1}-1}\br{\log\log x}^{\frac{n-1}{2n^{2}+1}+\varepsilon}}$, $n$ $\geq$ $4$, fully agrees with the best estimate in the $d$-dimensional, real hyperbolic manifolds with cusps case. Indeed, for $n$ $=$ $d$ and $\rho$ $=$ $\frac{1}{2}\br{n-1}$, the $O$-term in \cite[p. 3021, Th. 2]{AvdZenan} is given by
\begin{align*}
&O\br{x^{\alpha_{n}}\br{\log x}^{\beta_{n}-1}\br{\log\log x}^{\beta_{n}+\varepsilon}}\\
=&O\br{x^{\br{n-1}\br{1-\frac{2n+1}{4n^{2}+2}}}\br{\log x}^{\frac{n-1}{2n^{2}+1}-1}\br{\log\log x}^{\frac{n-1}{2n^{2}+1}+\varepsilon}}\\
=&O\br{x^{\br{n-1}-\frac{1}{2}\br{n-1}\frac{2n+1}{2n^{2}+1}}\br{\log x}^{\frac{n-1}{2n^{2}+1}-1}\br{\log\log x}^{\frac{n-1}{2n^{2}+1}+\varepsilon}}\\
=&O\br{x^{2\rho-\rho\frac{2n+1}{2n^{2}+1}}\br{\log x}^{\frac{n-1}{2n^{2}+1}-1}\br{\log\log x}^{\frac{n-1}{2n^{2}+1}+\varepsilon}}.
\end{align*}

As noted earlier, the $O$-terms in (\ref{dvadesetosam}), and hence the elements of the sequence $\set{2\rho-\rho\frac{2j+1}{2nj+1}}_{j}$, stem from the application of the Gallagher-Koyama techniques to the explicit formulas for the counting functions $\psi_{j}\br{x}$, $j$ $\geq$ $n$, given by Theorem \ref{Th1}. This fact, and the fact that the error terms are decreasing when $j$ is decreasing, motivated us to consider the possibility of further reduction of the error term in (\ref{dvadesetosam}) via the function $\psi_{n-1}\br{x}$. Unfortunately, we were not able to derive the explicit formula for $\psi_{n-1}\br{x}$ which would be analogous to (\ref{cetiri}) (see, remark at the end of Subsection 3.1). However, we were in position to adjust Hejhal's methods \cite{Hejhal1,Hejhal2}, to our setting, and thus obtain the explicit formula for $\psi_{n-1}\br{x}$, which is conditional and asymptotic, and is given by Theorem \ref{Th4}. In particular, we required the integral $\int\limits_{-A+\II{}T}^{\rho+\II{}T}-\frac{Z_{R}^{'}\br{s}}{Z_{R}\br{s}}\frac{x^{s+n-1}}{\prod\limits_{k=0}^{n-1}\br{s+k}}ds$ to be $O\br{\frac{x^{2\rho+\varepsilon_{1}+n-1}}{\varepsilon_{1}T^{1-\delta}}}$, where $\varepsilon_{1}$ $>$ $0$, $\delta$ $>$ $0$, and $T$ $\gg$ $0$, $A$ $\gg$ $2\rho$ are selected such that $\II{}T$ is not a singularity of $Z_{S}\br{s,\tau}$, $\tau$ $\in$ $\mathcal{T}$, and $-A$ is not a pole of $-\frac{Z_{R}^{'}\br{s}}{Z_{R}\br{s}}$. Note that an analogous condition of the given one is easily verified to be true in the case of compact Riemann surfaces for example ($n$ $=$ $2$, $\rho$ $=$ $\frac{1}{2}\br{n-1}$ $=$ $\frac{1}{2}$). Namely, according to Propositions 6.12-6.14 in \cite[pp. 105-108]{Hejhal1}, and the discussion given on page 108 of the same book, the integral $\int\limits_{-A+\II{}T}^{1+\varepsilon+\II{}T}\frac{Z^{'}\br{s}}{Z\br{s}}\frac{x^{s+1}}{s\br{s+1}}ds$ is $O\br{\frac{x^{2+\varepsilon}}{\varepsilon T}}$, where $Z\br{s}$ is the corresponding Selberg zeta function. The fact that the Ruelle zeta function in our case is a product of the Selberg zeta functions complicates things, so we left the aforementioned condition as an assumption.

As expected, the error term corresponding to $\psi_{n-1}\br{x}$ is\\ $O\br{x^{2\rho-\rho\frac{2\br{n-1}+1}{2n\br{n-1}+1}}\br{\log x}^{\frac{n-1}{2n\br{n-1}+1}-1}\br{\log\log x}^{\frac{n-1}{2n\br{n-1}+1}+\varepsilon}}$.

This is in line with (\ref{dvadesetosam}).

The obtained error term improves the error terms in (\ref{dvadesetosam}).

In particular, it improves the error term derived in \cite[p. 9, Th. 3]{MDPI1}.

If $n$ $=$ $2$, $\rho$ $=$ $\frac{1}{2}\br{n-1}$ $=$ $\frac{1}{2}$, our bound\\ $O\br{x^{2\rho-\rho\frac{2\br{n-1}+1}{2n\br{n-1}+1}}\br{\log x}^{\frac{n-1}{2n\br{n-1}+1}-1}\br{\log\log x}^{\frac{n-1}{2n\br{n-1}+1}+\varepsilon}}$ becomes\\ $O\br{x^{\frac{7}{10}}\br{\log x}^{\frac{1}{5}-1}\br{\log\log x}^{\frac{1}{5}+\varepsilon}}$, which is exactly the bound obtained by moving from the level $\psi_{0}\br{x}$ to the level $\pi_{\Gamma}\br{x}$ in the non-compact, hyperbolic, finite-volume Riemann surface with cusps case (Cf. \cite[p. 28, Th. 3.1]{AvdAnalMath}).

For $n$ $=$ $3$, $\rho$ $=$ $\frac{1}{2}\br{n-1}$ $=$ $1$,\\ $O\br{x^{2\rho-\rho\frac{2\br{n-1}+1}{2n\br{n-1}+1}}\br{\log x}^{\frac{n-1}{2n\br{n-1}+1}-1}\br{\log\log x}^{\frac{n-1}{2n\br{n-1}+1}+\varepsilon}}$ is\\ $O\br{x^{\frac{21}{13}}\br{\log x}^{-\frac{11}{13}}\br{\log\log x}^{\frac{2}{13}+\varepsilon}}$. This is in line with the result in the hyperbolic three manifolds case (see, \cite[p. 691, Th. 1.2]{Erata}).

In the remark on lower dimensions \cite[p. 3025]{AvdZenan}, the authors noted that one is in position to derive the asymptotics of $\psi_{0}\br{x}$ from $\psi_{d-1}\br{x}$ instead of reaching for $\psi_{d}\br{x}$ in the $d$-dimensional, real hyperbolic manifolds with cusps case. Taking into account the results we derived in this research, it becomes pretty obvious that the improved variant of Theorem 2 in \cite{AvdZenan}(the one obtained from $\psi_{d-1}\br{x}$), should state that for $\varepsilon$ $>$ $0$, there is a set $E$ of finite logarithmic measure, such that:
\begin{equation}\label{cetrdesetcetiri}
\begin{aligned}
&\pi_{\Gamma}\br{x}\\
=&\sum\limits_{\br{d-1}-\frac{1}{2}\br{d-1}\frac{2\br{d-1}+1}{2d\br{d-1}+1}<s_{n}\br{k}\leq2d_{0}}\br{-1}^{k}\li\br{x^{s_{n}\br{k}}}+\\
&O\br{x^{\br{d-1}-\frac{1}{2}\br{d-1}\frac{2\br{d-1}+1}{2d\br{d-1}+1}}\br{\log x}^{\frac{d-1}{2d\br{d-1}+1}-1}\br{\log\log x}^{\frac{d-1}{2d\br{d-1}+1}+\varepsilon}}
\end{aligned}
\end{equation}
as $x$ $\rightarrow$ $\infty$, $x$ $\in$ $\mathbb{R}$ $\backslash$ $E$, where $d_{0}$ $=$ $\frac{d-1}{2}$, $\br{s_{n}\br{k}-k}\br{2d_{0}-k-s_{n}\br{k}}$ is a small eigenvalue in $\SqBr{0,\frac{3}{4}d_{0}^{2}}$ of $\Delta_{k}$ on $\pi_{\sigma_{k},\lambda_{n}\br{k}}$ with $s_{n}\br{k}$ $=$ $d_{0}$ $+$ $\II{}\lambda_{n}\br{k}$ or $s_{n}\br{k}$ $=$ $d_{0}$ $-$ $\II{}\lambda_{n}\br{k}$ in $\spbr{\frac{3}{2}d_{0},2d_{0}}$, $\Delta_{k}$ is the Laplacian acting on the space of $k$-forms over $Y$, and $\pi_{\sigma_{k},\lambda_{n}\br{k}}$ is the principal series representation.

Notice that we used the notation $Y$ to temporarily denote the $d$-dimensional, real hyperbolic manifold with cusps in the particular case (\ref{cetrdesetcetiri}).

Also, we point out that the $O$-term $O\br{x^{2\rho-\frac{\rho}{n}}\br{\log x}^{-1}}$ in our unconditional prime geodesic theorem (\ref{cetrdesettri}) does not depend on the choice of the counting function $\psi_{j}\br{x}$, $j$ $\geq$ $n$ (Cf. \cite[p. 6, relation (8)]{MDPI1}).

\end{document}